%% file: definable.tex
\pdfoutput=1
\documentclass[12pt,a4paper,leqno]{article}

\usepackage[centertags]{amsmath}
\usepackage{amsthm}
\usepackage{amssymb,exscale}
\usepackage{paralist}
\usepackage{latexsym}
\usepackage{makeidx}
\usepackage[perpage,multiple,symbol]{footmisc}
\usepackage{xcolor}
\usepackage{hyperref}
\hypersetup{breaklinks=true}
\usepackage[backref]{enotez}
\DeclareTranslation{English}{enotez-title}{Notes and references}

\makeindex

\input{preamble}

\begin{document}
\title{Can each number be specified by a finite text?}
\author{Boris Tsirelson}
\date{}
\maketitle

\begin{abstract}
Contrary to popular misconception, the question in the title is far from simple. It involves sets of numbers on the
first level, sets of sets of numbers on the second level, and so on, endlessly. The infinite hierarchy of the levels
involved distinguishes the concept of ``definable number'' from such notions as ``natural number'', ``rational number'',
``algebraic number'', ``computable number'' etc.
\end{abstract}

\setcounter{tocdepth}{2}
\tableofcontents

\section{Introduction}
\label{sect1}

The question in the title may seem simple, but is able to cause controversy and trip up professional
mathematicians. Here is a quote from a talk ``Must there be numbers we cannot describe or define?'' \cite{Ha1}
by \href{http://en.wikipedia.org/wiki/Joel David Hamkins}{J.D. Hamkins}.\\[3pt]
\hbox{}\qquad\textcolor{blue!100}{The math tea argument}\\[3pt]
\hbox{}\qquad\quad Heard at a good math tea anywhere:\\[5pt]
\newlength\mylength
\setlength\mylength\textwidth
\addtolength\mylength{-4em}
\hbox{}\qquad\qquad\parbox{\mylength}{\noindent\emph{``There must be real numbers we cannot describe or define, because
    there are uncountably many real numbers, but only countably many definitions.''}}\\[5pt]
\hbox{}\qquad\quad Does this argument withstand scrutiny?\\[3pt]
See also ``Maybe there's no such thing as a random sequence'' \cite{Do} by \href{https://math.dartmouth.edu/~doyle/}{P.G. Doyle}
(in particular, on pages 6,7 note two excerpts from \href{http://en.wikipedia.org/wiki/Alfred Tarski}{A. Tarski}
\cite{Tar2}). And on Wikipedia one can also find the flawed ``math tea'' argument on talk pages and obsolete versions of
articles.\endnote{%
  \emph{Whichever definition of `definable' you choose, the formula that defines a definable object is a finite sequence of
  characters belonging to a finite alphabet. Thus the set of definable objects is definitively countable.} (From
  \href{https://en.wikipedia.org/w/index.php?Talk:Integer_sequence&oldid=824167431\#Definable_sequences}
  {Wikipedia}. A talk page. 2018.)}\,%
\endnote{%
  \emph{Classical mathematics permits (and requires) the existence of undefinable objects, but many people find this
  philosophically disquietening, questioning how an object can be said to exist if no mathematical statement can be used
  to uniquely identify it.\\As a result, a few mathematicians have developed systems of mathematics that do not
  involve undefinable objects.} (From
  \href{https://en.wikipedia.org/w/index.php?title=Definable&oldid=7518664}{Wikipedia}. Obsolete version of an article. 2004.)}
And elsewhere on the Internet.\endnote{%
  \emph{The describable numbers are all numbers for which there is any possible finite description that uniquely identifies
  the number. The countability argument still works: you can still enumerate all possible finite-length strings that
  could be descriptions, and define a one-to-one correspondence between strings that could be descriptions and natural
  numbers. The set of describable numbers is, thus, still countable, and the set of undescribables is not, which implies
  that the set of undescribables is far, far larger that the describables.}
  (From: \href{http://goodmath.scientopia.org/about-markcc/}{Chu-Carroll}, Mark C. (2014)
  \href{http://www.goodmath.org/blog/2014/05/26/you-cant-even-describe-most-numbers/}{``You can't even describe most
    numbers!''} ``Good Math/Bad Math'' blog.)}
I, the author, was myself a witness and accomplice. I shared and voiced the flawed
argument in informal discussions (but not articles or lectures). Despite some awareness (but not professionalism) in
mathematical logic,\endnote{%
  Tsirelson, Boris (2003)
  \href{https://www.tau.ac.il/~tsirel/Research/myspace/remins.html}{``Reminiscences''}. Self-published.}
I was a small part of the problem, and now I try to become a small part of the solution, spreading the truth.

Careless handling of the concept ``number specified by a finite text'' leads to paradoxes; in particular,
\href{http://en.wikipedia.org/wiki/Richard's paradox}{Richard's paradox}.\index{paradox, Richard's}\footnote{%
 \emph{The paradox begins with the observation that certain expressions of natural language define real numbers
  unambiguously, while other expressions of natural language do not. For example, ``The real number the integer part of
  which is 17 and the $n$th decimal place of which is 0 if $n$ is even and 1 if $n$ is odd'' defines the real
  number 17.1010101... = 1693/99, while the phrase ``the capital of England'' does not define a real number.\\ Thus there
  is an infinite list of English phrases (such that each phrase is of finite length, but lengths vary in the list) that
  define real numbers unambiguously. We first arrange this list of phrases by increasing length, then order all phrases
  of equal length \href{http://en.wikipedia.org/wiki/Lexicographical order}{lexicographically} (in dictionary order),
  so that the ordering is \href{http://en.wikipedia.org/wiki/Canonical form}{canonical}. This yields an infinite list
  of the corresponding real numbers: $r_1,r_2,\dots$  Now define a new real number $r$ as follows. The integer part of
  $r$ is $0$, the $n$th decimal place of $r$ is $1$ if the $n$th decimal place of $r_n$ is not $1$, and the $n$th decimal
  place of $r$ is $2$ if the $n$th decimal place of $r_n$ is 1.\\ The preceding two paragraphs are an expression in English
  that unambiguously defines a real number $r$. Thus $r$ must be one of the numbers $r_n$. However, $r$ was constructed
  so that it cannot equal any of the $r_n$. This is the paradoxical contradiction.} (Quoted from Wikipedia.)}
See also
\href{https://www.dpmms.cam.ac.uk/~wtg10/richardsparadox.html}{``Definability paradoxes''} by Timothy Gowers.

In order to ask (and hopefully solve) a well-posed question we have to formalize the concept ``number specified by a
finite text'' via a well-defined mathematical notion ``definable number''. What exactly is meant by ``text''? And what
exactly is meant by ``number specified by text''? Does ``specified'' mean ``defined''? Can we define such notions as
``definition'' and ``definable''? Striving to understand definitions in general, let us start with some examples.

136 notable constants are collected, defined and discussed in the book ``Mathematical constants'' by Steven
Finch \cite{Fi}.
The first member of this collection is \href{http://en.wikipedia.org/wiki/Square root
  of 2}{``Pythagoras' Constant, $\sqrt2$''};
the second is \href{http://en.wikipedia.org/wiki/Golden ratio}{``The Golden Mean, $\varphi$''}; the third
\href{http://en.wikipedia.org/wiki/E (mathematical constant)}{``The Natural Logarithmic Base, $e$''}; the fourth
\href{http://en.wikipedia.org/wiki/Pi}{``Archimedes' Constant, $\pi$''}; and the last (eleventh) in Chapter 1
``Well-Known Constants'' is \href{http://en.wikipedia.org/wiki/Chaitin's constant}{``Chaitin's Constant''}.

Each constant has several equivalent definitions. Below we take for each constant the first (main) definition from the
mentioned book.
\begin{compactitem}
\item The first constant $\sqrt2$ is defined as the positive real number whose product by itself is equal to 2. That
is, the real number $x$ satisfying $x>0$ and $x^2=2$.
\item The second constant $\varphi$ is defined as the real number satisfying $\varphi>0$ and
$1+\frac1\varphi=\varphi$.
\item The third constant $e$ is defined as the limit of $\textstyle (1+x)^{1/x}$ as $x\to0$. That is, the real
number satisfying the following condition:
\end{compactitem}
\begin{quote}for every $\varepsilon>0$ there exists $\delta>0$ such that for every $x$ satisfying $-\delta<x<\delta$ and
  $x\ne0$ holds $\textstyle -\varepsilon<(1+x)^{1/x}-e<\varepsilon$.\end{quote}
The same condition in symbols:\footnote{%
  Logical notation: \hfill $\land$ \, "and" \hfill $\lor$ \, "or" \hfill $\imp$ \, "implies" \hfill $\neg$ \,
  "not" \linebreak
  $\forall$ \, "for every" \hfill $\exists$ \, "there exists (at least one)" \hfill $\exists!$ \, "there exists one and
  only one" \linebreak
  \href{http://en.wikipedia.org/wiki/List of logic symbols}{(link to a longer list)}.}
\begin{multline*}
\forall \varepsilon>0 \;\> \exists \delta>0 \;\> \forall x \;\; \(\, ( -\delta<x<\delta \,\,\land\,\, x\ne0 )
 \imp \\
( -\varepsilon<(1+x)^{1/x}-e<\varepsilon ) \,\).
\end{multline*}

We note that these three definitions are of the form ``the real number $x$ satisfying $P(x)$'' where $P(x)$ is a statement
that may be true or false depending on the value of its variable $x$; in other words, not a statement when $x$ is just a
variable, but a statement whenever a real number is substituted for the variable. Such $P(x)$ is called a property of
$x$, or a \href{http://en.wikipedia.org/wiki/Predicate (mathematical logic)}{predicate} (on real
numbers).\index{predicate}

Not all predicates may be used this way. For example, we cannot say ``the real number $x$ satisfying $x^2=2$'' (why
``the''? two numbers satisfy, one positive, one negative), nor ``the real number $x$ satisfying $x^2=-2$'' (no such
numbers). In order to say ``the real number $x$ satisfying $P(x)$'' we have to prove existence and uniqueness:

\begin{quote}
existence: $\exists x \;\> P(x)$ (in words: there exists $x$ such that $P(x)$);

uniqueness: $\forall x,y \;\> \(\, ( P(x) \land P(y) ) \imp (x=y) \,\)$\\
\hbox{}\quad(in words: whenever $x$ and $y$ satisfy $P$ they are equal).
\end{quote}

In this case one
\href{http://en.wikipedia.org/wiki/Uniqueness quantification}{says ``there is one and only one such $x$'' and writes
``$\exists! x \; P(x)$''}.

The road to definable numbers passes through definable predicates. We postpone this matter to the next section and
return to examples.
\begin{compactitem}
\item The fourth constant $\pi$ is defined as the area enclosed by a circle of radius 1.
\end{compactitem}
This definition involves geometry. True, a lot of equivalent definitions in terms of numbers are well-known; in
particular, according to the mentioned book, this area is equal to $\textstyle \, 4\int_0^1 \sqrt{1-x^2} \, dx =
\lim_{n\to\infty} \frac4{n^2} \sum_{k=0}^n \sqrt{n^2-k^2} \, $. However, in general, every branch of mathematics may be
involved in a definition of a number; existence of an equivalent definition in terms of (only) numbers is not
guaranteed.

The last example is Chaitin's constant.\index{Chaitin's constant} In contrast to the four constants (mentioned above) of evident
theoretical and practical importance, Chaitin's constant is rather of theoretical interest.  Its definition is
intricate. Here is a simplified version, sufficient for our purpose.\endnote{%
 % Not quite the Chaitin's constant?
 \emph{Because $\Om$ depends on the program encoding used, it is sometimes called Chaitin's construction instead of
   Chaitin's constant when not referring to any specific encoding.} (From
 Wikipedia:\href{http://en.wikipedia.org/wiki/Chaitin's  constant}{Chaitin's Constant}.)
 The polynomial $f$ used in this definition is not uniquely determined as ``such that the sequence $A_1, A_2, \dots$ is
 uncomputable'', but for now every such polynomial serves our purpose; later, in Section \ref{sect8}, we'll add one more
 requirement related to the uncomputability. Other requirements, related to randomness, are irrelevant to this
 essay. Existence of such $f$ follows from Matiyasevich's theorem; 9 unknowns are sufficient according to: James Jones
 (1982) \href{https://doi.org/10.2307/2273588}{``Universal Diophantine equation''} \emph{J.~Symbolic Logic} (Cambridge)
 \textbf{47} (3): 549--571.}
 
\begin{sloppypar}
\begin{compactitem}
\item The last constant $\Om$ is defined as the sum of the series $\textstyle \Om = \sum_{N=1}^\infty 2^{-N}
A_N$ where $A_N$ is equal to 1 if there exist natural numbers $x_1,x_2,x_3,x_4,x_5,x_6,x_7,x_8,x_9$ such that
$f(N,x_1,x_2,x_3,x_4,x_5,x_6,x_7,x_8,x_9)=0$, otherwise $A_N=0$; and $f$ is a polynomial in 10 variables, with integer
coefficients, such that the sequence $A_1, A_2, \dots$ is uncomputable.
\end{compactitem}
\href{http://en.wikipedia.org/wiki/Hilbert’s tenth problem}{Hilbert's tenth problem} asked for a general algorithm
that could ascertain whether the \href{http://en.wikipedia.org/wiki/Diophantine equation}{Diophantine equation}
$f(x_0,\dots,x_k)=0$ has positive integer solutions $(x_0,\dots,x_k)$, given arbitrary polynomial $f$ with integer
coefficients. It appears that no such algorithm can exist even for a single $f$ and arbitrary $x_0$, when $f$ is
complicated enough. See Wikipedia: \href{http://en.wikipedia.org/wiki/Computability theory}{computability theory},
\href{http://en.wikipedia.org/wiki/Diophantine set\#Matiyasevich's_theorem}{Matiyasevich's theorem}; and
Scholarpedia:\href{http://www.scholarpedia.org/article/Matiyasevich_theorem}{Matiyasevich theorem}.\index{Matiyasevich theorem}
\end{sloppypar}

The five numbers $\sqrt2, \varphi, e, \pi, \Om$ are defined, thus, should be definable according to any reasonable
approach to definability. The first four numbers $\sqrt2, \varphi, e, \pi$ are computable\index{computable!number}\index{number!computable} (both theoretically and
practically; in fact, trillions, that is, millions of millions, of decimal digits of $\pi$ are already computed), but
the last \index{number!uncomputable}number $\Om$ is uncomputable. How so? Striving to better understand this strange situation we may introduce
approximations $A_{M,N}$ to the numbers $A_N$ as follows: $A_{M,N}$ is equal to 1 if there exist natural numbers
$x_1,x_2,x_3,x_4,x_5,x_6,x_7,x_8,x_9$ less than $M$ such that $f(N,x_1,x_2,x_3,x_4,x_5,x_6,x_7,x_8,x_9)=0$, otherwise
$A_N=0$; here $M$ is arbitrary. For each $N$ we have $A_{M,N}\uparrow A_N$ as $M\to\infty$; that is, the sequence
$A_{1,N}, A_{2,N}, \dots$ is increasing, and converges to $A_N$. Also, this sequence $A_{1,N}, A_{2,N}, \dots$ is
computable (given $M$, just check all the $(M-1)^9$ points $(N,x_1,x_2,x_3,x_4,x_5,x_6,x_7,x_8,x_9)$, $0<x_1<M, \dots,
0<x_9<M$). Now we introduce approximations $\om_M$ to the number $\Om$ as follows: $\textstyle \om_M = \sum_{N=1}^M
2^{-N} A_{M,N}$. We have $\om_M\uparrow\Om$ (as $M\to\infty$), and the sequence $\om_1, \om_2, \dots$ is computable. A
wonder: a computable increasing sequence of rational numbers converges to a uncomputable number!

For every $N$ there exists $M$ such that $A_{M,N}=A_N$; such $M$ depending on $N$, denote it $M_N$ and get $\textstyle
\sum_{N=1}^\infty 2^{-N} A_{M_N,N}=\Om$; moreover, $\textstyle \Om-\sum_{N=1}^K 2^{-N} A_{M_N,N}\le 2^{-K}$ for all
$K$. In order to compute $\Om$ up to $2^{-K}$ it suffices to compute $\textstyle \sum_{N=1}^K 2^{-N} A_{M_N,N}$.
Doesn't it mean that $\Om$ is computable? No, it does not, unless the sequence $M_1,M_2,\dots$ is computable. Well,
these numbers need not be optimal, just large enough. Isn't $\textstyle M_N=10^{1000N}$ large enough? Amazingly, no,
this is not large enough. Moreover, $\textstyle M_N=10^{10^{1000N}}$ is not enough. And even the ``power tower''
$M_N=\underbrace{10^{10^{\cdot^{\cdot^{10}}}}}_{1000N}$ is still not enough!

Here is the first paragraph from a prize-winning article by Bjorn Poonen \cite{Po}:
\begin{quote}
\emph{Does the equation $\textstyle x^3+y^3+z^3=29$ have a solution in integers? Yes: $(3, 1,
1)$, for instance. How about $\textstyle x^3+y^3+z^3=30$? Again yes, although this was not known until 1999: the smallest
solution is $(-283059965, -2218888517, 2220422932)$. And how about $\textstyle x^3+y^3+z^3=33?$ This is an unsolved
problem.}
\end{quote}
Given that the simple Diophantine equation $\textstyle N+x^3+y^3-z^3=0$ has solutions for $N=30$ but only
beyond $10^9$ we may guess that the "worst case" Diophantine equation $f(N,x_1,x_2,x_3,x_4,x_5,x_6,x_7,x_8,x_9)=0$ needs
\emph{very} large $M_N$. In fact, the sequence $M_1,M_2,\dots$ has to be uncomputable (otherwise $\Om$ would be
computable, but it is not). Some computable sequences grow fantastically fast. See Wikipedia:
\href{http://en.wikipedia.org/wiki/Ackermann function}{``Ackermann function''},
\href{http://en.wikipedia.org/wiki/Fast-growing hierarchy}{``Fast-growing hierarchy''}.  And nevertheless, no one of them
bounds from above the sequence $M_1,M_2,\dots\,$ Reality beyond imagination!

Every computable number is definable, but a definable number need not be computable. Computability being another story,
we return to definability.

\section{From predicates to relations}
\label{sect2}

Recall the five definitions mentioned in the introduction. They should be special cases of a general notion
``definition''. In order to formalize this idea we have to be more pedantic than in the introduction. \emph{``Nothing but the hard technical story is any real good''} (\href{http://en.wikipedia.org/wiki/John Edensor Littlewood}{Littlewood}, \href{http://en.wikipedia.org/wiki/A Mathematician's Miscellany}{A Mathematician's Miscellany}, page 70); exercises are waiting for you.

All mathematical objects (real numbers, limits, sets etc.) are treated in the framework of the mainstream mathematics,
unless stated otherwise. Alternative approaches are sometimes mentioned in Sections \ref{sect9},
\ref{sect10}. \href{http://en.wikipedia.org/wiki/Naive_set_theory}{Naive set theory} suffices for Sections
\ref{sect2}--\ref{sect7}; \href{http://en.wikipedia.org/wiki/Set_theory#Axiomatic_set_theory}{axiomatic set theory} is used in
Sections \ref{sect8}--\ref{sect10}.

A definition is a text in a language. A straightforward formalization of such notions as ``definition'' and ``definable''
uses ``formal language'' (a formalization of ``language'') and other notions of \href{http://en.wikipedia.org/wiki/Model
  theory}{model theory}. Surprisingly, there is a shorter
way. \href{http://en.wikipedia.org/wiki/Set_(mathematics)\#Basic_operations}{Operations on sets} are used instead of
\href{http://en.wikipedia.org/wiki/First-order_logic\#Logical_symbols}{logical symbols}, and
\href{http://en.wikipedia.org/wiki/Finitary relation}{relations} instead of
\href{http://en.wikipedia.org/wiki/Predicate (mathematical logic)}{predicates}.

\emph{``However, predicates have many different uses and interpretations in mathematics and logic, and their precise
definition, meaning and use will vary from theory to theory.''} (\href{http://en.wikipedia.org/wiki/Predicate
  (mathematical logic)}{Quoted from Wikipedia.}) Here we use predicates for informal explanations only; on the formal
level they will be avoided (replaced with relations).

The number $\sqrt2$ was defined as the real number $x$ such that $P(x)$, where $P(x)$ is the predicate ``$x>0$ and
$x^2=2$''. This predicate is the conjunction $P_1(x)\land P_2(x)$ of two predicates $P_1(x)$ and $P_2(x)$, the first
being ``$x>0$'', the second ``$x^2=2$''. The single-element set $A=\{x\in\R\mid P(x)\}=\{\sqrt2\}$ corresponding to
the predicate $P(x)$ is the intersection $A=A_1\cap A_2$ of the sets $A_1=\{x\in\R\mid P_1(x)\}=(0,\infty)$ and
$A_2=\{x\in\R\mid P_2(x)\}=\{-\sqrt2,\sqrt2\}$. (Here and everywhere, $\R$ is the set of all real
numbers.)\index{zzR@$\R$, the set of all real numbers}%
\footnote{%
  Set notation: \\
  $A=\{x\mid P(x)\}$ \, \href{http://en.wikipedia.org/wiki/Set-builder notation\#Sets_defined_by_a_predicate}{``$A$ is
    the set of all $x$ such that $P(x)$''} \hfill
  $x\in A$ \, ``$x$ \href{http://en.wikipedia.org/wiki/Element (mathematics)}{belongs} to $A$'' \\
  $A\cup B$ \, \href{http://en.wikipedia.org/wiki/Union (set theory)}{union} \hfill
  $A\cap B$ \, \href{http://en.wikipedia.org/wiki/Intersection (set theory)}{intersection} \hfill
  $A\setminus B$ \, \href{http://en.wikipedia.org/wiki/Complement (set theory)\#Relative_complement}{set difference}
   \hfill
  $A\times B$ \, \href{http://en.wikipedia.org/wiki/Cartesian product}{Cartesian product} \\
  $\R$ \, \href{http://en.wikipedia.org/wiki/Real line}{real line} \hfill
  $\R^2$ \, \href{http://en.wikipedia.org/wiki/Real coordinate space}{Cartesian plane} \hfill\hfill
  and \href{http://en.wikipedia.org/wiki/Set (mathematics)}{more}, \,
  \href{http://en.wikipedia.org/wiki/Algebra_of_sets}{more}, \,
  \href{http://en.wikipedia.org/wiki/Set-builder_notation}{more}.\hfill\hbox{}}

This is instructive. \emph{In order to formalize a definition of a number via its defining property, we have to deal with
sets of numbers, and more generally, relations between numbers.}

Also, $x^2$ is the product $x\cdot x$, and $2$ is the sum $1+1$. But what is ``product'', ``sum'', ``1'' and ``0''? The
answer is given by the \href{http://en.wikipedia.org/wiki/Real number\#Axiomatic_approach}{axiomatic approach to real
  numbers}: they are a complete totally ordered field.  It means that addition, multiplication and order are defined and
have the appropriate properties. Thus, 0 is defined as the real number $x$ satisfying the condition $\forall y \;
(x+y=y)$. Similarly, 1 is defined as the real number $x$ satisfying the condition $\forall y \; (x\cdot y=y)$.

Now we need predicates with two and more variables. The order is a binary\index{predicate!binary} (that is, with two variables) predicate ``$x\le
y$''. Addition is a \index{predicate!ternary}ternary (that is, with three variables) predicate ``$x+y=z$''. Similarly, multiplication is a ternary
predicate ``$xy=z$'' (denoted also ``$x\cdot y=z$'' or ``$x\times y=z$'').

Each unary\index{predicate!unary} (that is, with one variable) predicate $P(x)$ on real numbers leads to a set $\{x\in\R\mid P(x)\}$ of
real numbers, a subset of the real line $\R$. Likewise, each binary predicate $P(x,y)$ on reals leads to a set
$\{(x,y)\in\R^2\mid P(x,y)\}$ of pairs of real numbers, a subset of the Cartesian plane $\R^2=\R
\times \R$, the latter being the \href{http://en.wikipedia.org/wiki/Cartesian product\#A_two-dimensional_coordinate_system}{Cartesian product}\index{product of sets, Cartesian} of the real line by itself. On the other hand, a
\href{http://en.wikipedia.org/wiki/Binary relation}{binary relation}\index{relation!binary} on $\R$ is defined as an arbitrary subset
of $\R^2$.

Thus, each binary predicate on reals leads to a binary relation on reals. If we swap the variables, that is, turn to
another predicate $Q(x,y)$ that is $P(y,x)$, then we get another relation $\{(x,y)\mid Q(x,y)\}=\{(x,y)\mid
P(y,x)\}=\{(y,x)\mid P(x,y)\}$, \href{http://en.wikipedia.org/wiki/Inverse relation}{inverse} (in other words, converse,
or opposite) to the former relation (generally different, but sometimes the same).

Similarly, each \index{relation!ternary}ternary predicate on reals leads to a ternary relation on reals; and, changing the order of variables,
we get $3!=6$ ternary relations (generally, different) corresponding to 6
\href{http://en.wikipedia.org/wiki/Permutation}{permutations} of 3 variables. And generally, each \ary{n}
predicate on reals leads to a \ary{n} relation\index{relation!n@$n$-ary} on reals (a subset of $\R^n$); and, changing the order
of variables, we get $n!$ such relations. The case $n=1$ is included (for unification); a unary relation\index{relation!unary} on reals
(called also property of reals) is a subset of $\R$.

Thus, on reals, the order\index{relation!order} is the binary relation $\{(x,y)\mid x\le y\}$, the addition\index{relation!addition} is the ternary
relation$\{(x,y,z)\mid x+y=z\}$, and the multiplication\index{relation!multiplication} is the ternary relation $\{(x,y,z)\mid xy=z\}$. Still, we cannot
forget predicates until we understand how to construct new relations out of these basic relations. For example, how to
construct the binary relation $\{(x,y)\mid x+y=y\}$ and the unary relation $\{x\mid \forall y \; (x+y=y)\}\,$? We know
that if a predicate $P(x)$ is the conjunction $P_1(x)\land P_2(x)$ of two predicates, then it leads to the intersection
$A=A_1\cap A_2$ of the corresponding sets. Similarly, the disjunction $P_1(x)\lor P_2(x)$ leads to the union $A=A_1\cup
A_2$, and the negation $\neg P_1(x)$ leads to the complement $A=\R\setminus A_1$. Also, the implication $P_1(x)\imp
P_2(x)$ leads to $A=(\mathbb R\setminus A_1)\cup A_2$, and the equivalence $P_1(x)\equ P_2(x)$ leads to
$A=\((\R\setminus A_1)\cap(\R\setminus A_2)\) \cup (A_1\cap A_2)$. The same holds for \ary{n} predicates;
the disjunction $P_1(x,y)\lor P_2(x,y)$ still corresponds to the union $A=A_1\cup A_2$, the negation $\neg P_1(x,y,z)$
to the complement $A=\R^3\setminus A_1$, etc. But what to do when $P(x,y)$ is $P_1(x,y,y)$, or $P(x)$ is $\forall y \;
P_1(x,y)$, or $P(x,y,z)$ is $P_1(y,x)\land P_2(y,z)$, etc?

This question was answered, in context of axiomatic set theory, in the first half of the 20th century.\endnote{%
  See Wikipedia:\href{http://en.wikipedia.org/wiki/Von Neumann–Bernays–Godel set theory\#Class_existence_axioms_and_axiom_of_regularity}{Von Neumann-Bernays-G\"odel set theory}
 % (Section ``Class existence axioms and axiom of regularity'')
  and \href{http://en.wikipedia.org/wiki/Godel operation}{G\"odel operation}.}
A somewhat different answer, in the context of definability, was given by van den Dries in 1998 \cite{Dr1},\cite{Dr2} and
slightly modified by Auke Bart Booij in 2013 \cite{Bo}; see also Macintyre 2016 \cite[``Defining
First-Order Definability'']{Mac}. Here is the answer (slightly modified).

First, in addition to the Boolean operations\index{operation!Boolean (union, complement)} (union and complement; intersection is superfluous, since it is complement
of the union of complements) on subsets of $\R^n$, we introduce \index{operation!permutation}permutation of coordinates; for example ($n=3$),
$A=\{(x,y,z)\in\R^3\mid(z,x,y)\in A_1\}$; and in general,
\begin{equation*}
A=\{(x_1,\dots,x_n)\in\R^n\mid(x_{i_1},\dots,x_{i_n})\in A_1\}
\end{equation*}
where $(i_1,\dots,i_n)$ is an arbitrary permutation of $(1,\dots,n)$.

In particular, permutation of coordinates in a binary relation gives the inverse relation. For example, the inverse to
$\{(x,y)\mid x\le y\}$ is $\{(x,y)\mid y\le x\}=\{(x,y)\mid x\ge y\}$. And, by the way, the intersection of these two is
the relation $\{(x,y)\mid x=y\}$ (corresponding to the predicate ``$x=y$'').

Second, set multiplication,\index{operation!set multiplication} in other words, \href{http://en.wikipedia.org/wiki/Cartesian product}{Cartesian product}
by $\R$: $A=A_1\times\R$, that is,
\begin{equation*}
A=\{(x_1,\dots,x_{n+1})\in\R^{n+1}\mid(x_1,\dots,x_n)\in A_1\},
\end{equation*}
turns a \ary{n} relation to a relation that is formally \ary{(n+1)}, but the last variable is
unrelated to others.

Now, returning to a predicate $P(x,y,z)$ of the form $P_1(y,x)\land P_2(y,z)$, we treat the corresponding ternary
relation $A=\{(x,y,z)\in\R^3\mid P(x,y,z)\}$ as the intersection of two ternary relations $A_1=\{(x,y,z)\in\R^3\mid
P_1(y,x)\}$ and $A_2=\{(x,y,z)\in\R^3\mid P_2(y,z)\}$; and $A_1$ as the Cartesian product of the binary relation
$B_1=\{(x,y)\in\R^2\mid P_1(y,x)\}$ by $\R$, $B_1$ being inverse to the relation $B_2=\{(x,y)\in\R^2\mid P_1(x,y)\}$
(corresponding to the given predicate $P_1(x,y)$); and $A_2$ as obtained (by permutation of coordinates) from the
Cartesian product $\{(y,z,x)\in\R^3\mid P_2(y,z)\} = \{(y,z)\in\R^2\mid P_2(y,z)\}\times\R$ (by $\R$) of the relation
corresponding to the given predicate $P_2(y,z)$.

Third, the projection;\index{operation!projection} for example ($n=1$), $A=\{x\mid\exists y\in\R \; \((x,y)\in A_1\)\}$; and in
general,
\begin{equation*}
A=\{(x_1,\dots,x_n)\mid\exists x_{n+1}\in\R \; \((x_1,\dots,x_{n+1})\in A_1\)\};
\end{equation*}
it turns a $(n+1)$-ary relation to a $n$-ary relation. For $n=1$ the set $A$ is also called the
\href{http://en.wikipedia.org/wiki/Binary relation\#Formal_definition}{domain of the binary relation} $A_1$.

Now, returning to a predicate $P(x,y)$ of the form $P_1(x,y,y),\,$ we rewrite it as ``$\exists z\; \(P_1(x,y,z) \land
y=z\)$'' and treat the corresponding binary relation as the projection of the ternary relation $\{(x,y,z)\mid
P_1(x,y,z)\} \cap \{(x,y,z)\mid y=z\}$, and $\{(x,y,z)\mid y=z\}$ as a permutation of the Cartesian product $\{(y,z)\mid
y=z\}\times\R$.

What if $P(x)$ is ``$\forall y \; P_1(x,y)$''? Then we rewrite it as ``$\neg\exists y \; \neg P_1(x,y)$'' and get the
complement of the projection of the complement of the relation corresponding to $P_1(x,y)$.

So, we accept the 3 given relations (order, addition, multiplication) as ``definable'', and we accept the 5 operations
(complement, union, permutation, set multiplication, projection) for producing definable relations\index{relation!definable}\index{definable!relation} out of other
definable relations. Thus we get infinitely many definable relations (unary, binary, ternary and so on).

More formally, these relations are called ``first-order definable (without parameters) over $(\R; \le, +, \times)$'';
but, less formally, ``definable over'' is often replaced with ``definable in'' (and sometimes ``definable from'');
``without parameters'' is omitted throughout this essay; also ``first order'' and ``over $(\R; \dots)$'' are often
omitted in this section. See Wikipedia: ``Definable set'': \href{http://en.wikipedia.org/wiki/Definable
 set\#Definition}{Definition}; \href{http://en.wikipedia.org/wiki/Definable set\#The_field_of_real_numbers}{The field of
 real numbers}.

Generally, starting from a set (not necessarily the real line) and some chosen relations on this set (including the equality relation if needed), and applying the 5
operations (complement, union, permutation, set multiplication, projection) repeatedly (in all possible combinations),
one obtains an infinite collection of relations (unary, binary, ternary and so on) on the given set. Every such
collection of relations is called a structure (Booij \cite{Bo}), or a VDD-structure
(\href{https://www.maths.tcd.ie/~btyrrel/}{Brian Tyrrell} \cite{Ty}) on the given set. According to Tyrrell \cite[page
  3]{Ty}, ``The advantage of this definition is that no model theory is then needed to develop the theory''.  The
technical term ``VDD structure'' (rather than just ``structure'' used by van den Dries and Booij) is chosen by Tyrrell ``to
prevent a notation clash'' (Tyrrell \cite[page 2]{Ty}), since many other structures of different kinds are widely used in
mathematics. ``VDD'' apparently refers to van den Dries who pioneered this approach. But let us take a shorter term
``\Dstructure'',\index{D-structure (generated by)} where ``D'' refers to ``definable'' and ``Dries'' as well. The \Dstructure\ obtained (by the
5 operations) from the chosen relations, in other words, generated by these relations, is the smallest
\Dstructure\ containing these relations.

Generality aside, we return to the special case, the \Dstructure\ of definable relations on the real line defined above
(generated by order, addition and multiplication; though, the order appears to be superfluous).

\begin{exercise}
 Prove that a relation is definable in $(\R; \le, +, \times)$ if and only if it is definable in $(\R;
 +, \times)$. \emph{Hint:} $x\le y$ if and only if $\exists z\in\R \; (x+z^2=y)$.
\end{exercise}

We say that a number $x$ is definable,\index{definable!number}\index{number!definable} if the single-element set $\{x\}$ is a definable unary relation.
\begin{exercise}
 Prove that the numbers 0 and 1 are definable. \emph{Hint:} recall ``$\forall y \; (x+y=y)$'' and ``$\forall y \;
 (x\cdot y=y)$''.
\end{exercise}

\begin{exercise}
 Prove that the sum of two definable numbers is definable. \emph{Hint:} $\exists y\in\R \; \exists z\in\R
 \;\,  \((y\in A_1) \land (z\in A_2) \land (y+z=x)\)$.
\end{exercise}

\begin{exercise}
 Prove that the number $\frac{355}{113}$ is definable. \emph{Hint:} $\exists y\in\R \; \exists z\in\R \;\,
 (y=113 \land z=355 \land xy=z)$.
\end{exercise}

\begin{exercise}
 Prove that the number $\sqrt2$ is definable. \emph{Hint:} $(x>0)\land(x\cdot x=2)$.
\end{exercise}

\begin{exercise}
 Prove that the golden ratio $\varphi$ is definable.
\end{exercise}

\begin{exercise}
 Prove that the binary relation ``$y=|x|$'' is definable. \emph{Hint:} $(x^2=y^2 \land y\ge0)$.
\end{exercise}

In contrast, the ternary relation ``$x^y=z$'' is not definable.\index{relation!undefinable}\index{undefinable!relation} Moreover, the binary relation $\{(x,y)\mid y=2^x \land
0\le x\le1 \}$ is not definable. The problem is that all relations definable in $(\R; +, \times)$ are
\href{http://en.wikipedia.org/wiki/Semialgebraic set}{semialgebraic sets}\index{semialgebraic} over (the subring of) integers.\footnote{%
  \emph{Theorem.} All relations definable in $(\R; +, \times)$ are semialgebraic sets over integers.

  \emph{Proof.} The two relations ``$+$'', ``$\times$'' are semialgebraic (evidently). Two operations, permutation and
  set multiplication, applied to semialgebraic relations, give semialgebraic relations (evidently). The third,
  projection operation, applied to a semialgebraic relation, gives a semialgebraic relation by the
  \href{http://en.wikipedia.org/wiki/Tarski–Seidenberg theorem}{Tarski-Seidenberg theorem}\cite[Theorem 2.76]{Ba}.

  \emph{Theorem.} If $a>1$ and $-\infty<b<c<\infty$, then the binary relation $\{(x,y)\mid (y=a^x) \land (b\le x\le c)
  \}$ is not semialgebraic.

  \emph{Proof.} Assume the contrary. Then the function $x\mapsto a^x$ on $[b,c]$, being semialgebraic, must be
  algebraic.\cite[Prop. 2.86]{Ba}, \cite[Corollary 3.5]{Mar}. It means existence of a polynomial $p(\cdot,\cdot)$ (not
  identically 0) such that $p(x,a^x)=0$ for all $x\in[b,c]$. It follows that $p(z,e^{z\log a})=0$ for all complex
  numbers $z$. Taking $z=\tfrac{2n\pi i}{\log a}$ we get $p\(\tfrac{2n\pi i}{\log a},1\)=0$ for all integer
  $n$. Therefore $p(z,1)=0$ for all complex $z$ (otherwise the polynomial $z\mapsto p(z,1)$ cannot have infinitely many
  roots). Similarly, taking $z=\tfrac{\log u+2n\pi i}{\log a}$ we get $p(z,u)=0$ for all complex $z$ and all $u>0$,
  therefore everywhere; a contradiction.}

Thus, we cannot define the number $e$ via $(1+x)^{1/x}$ in this framework.\index{number!undefinable}\index{undefinable!number}  Also, only
\href{http://en.wikipedia.org/wiki/Algebraic number}{algebraic numbers} are definable in this framework.

Each natural number is definable, which does not mean that the set $\N$\index{zzN@$\N$, the set of all natural numbers}
of all natural numbers is definable (in $(\R; +, \times)$). In fact, it is not!\footnote{%
  Follows immediately from the lemma below.

  \emph{Lemma.} For every semialgebraic subset $A$ of $\R$ there exists
  $a\in\R$ such that either $(a,\infty)\subset A$ or $(a,\infty)\cap A=\emptyset$.

  \emph{Proof.} First, the claim holds for every set of the form $A=\{x\in\R\mid p(x)>0\}$ where $p(\cdot)$ is a
  polynomial, since either $p(x)\to+\infty$ as $x\to+\infty$, or $p(x)\to-\infty$ as $x\to+\infty$, or $p(x)$ is
  constant. Second, a Boolean operation (union, complement), applied to sets that satisfy the claim, gives a set that
  satisfies the claim (evidently).}

\enlargethispage{3pt}

We could accept the set $\N$ of natural numbers as definable, that is, turn to definability in $(\R; +,\times,\N)$, but
does it help to define the number $e$? Surprisingly, it does! ``[\dots] then the situation changes drastically'' (van den
Dries \cite[Example 1.3]{Dr1}). See also Booij \cite[page 17]{Bo}: ``[\dots] if we add the seemingly innocent set $\mathbf
Z$ to the tame structure of semialgebraic sets, we get a wild structure [\dots]''

\section{Beyond the algebraic}
\label{sect3}

\begin{sloppypar}
\textbf{In this section, ``definable'' means ``first order definable in $(\R; +,\times,\N)$''.} In other words, the real
line is endowed with the \Dstructure\ generated by addition, multiplication, and the set of natural numbers. Good news:
we'll see that the five numbers $ \sqrt {2}, \varphi, e, \pi, \Om, $ discussed in Introduction, are definable. Bad news:
in addition to their usual definitions we'll use Diophantine equations, computability\index{computability} and Matiyasevich's theorem
(mentioned in Introduction in relation to Chaitin's constant). The reader not acquainted with computability theory
should rely on intuitive idea of computation (instead of formal proofs of computability), and consult the linked
Wikipedia article for computability-related notions (``recursively enumerable'', ``computable sequence''). Alternatively,
the reader may skip to Section \ref{sect5}; there, usual definitions will apply, no computability needed.
\end{sloppypar}

Every \index{set!Diophantine}\href{http://en.wikipedia.org/wiki/Diophantine set}{Diophantine set}
\begin{equation*}
\{(a_1,\dots,a_n)\in\N^n\mid \exists x_1,\dots,x_m\in\N \;\; p(a_1,\dots,a_n,\,x_1,\dots,x_m)=0 \}
\end{equation*}
(where $p(\dots)$ is a polynomial with integer coefficients), treated as a subset of $\R^n$, is a definable
\ary{n} relation. And every \index{set!recursively enumerable}\href{http://en.wikipedia.org/wiki/Recursively enumerable set}{recursively enumerable}
set is Diophantine.\index{set!Diophantine}

For every \href{http://en.wikipedia.org/wiki/Integer sequence\#Computable_and_definable_sequences}{computable sequence}\index{computable!sequence}
$(k_1,k_2,\dots)$ of natural numbers, the binary relation $\{(n,k)\in\N^2 \mid k=k_n\} = \{ (1,k_1), (2,k_2), \dots \}$
is recursively enumerable, therefore definable.

In particular, the binary relation ``$x\in\N$ and $y=x^x$'' is definable, as well as ``$x\in\N$ and $y=(x+1)^x$''. Now
(at last!) the number $e$ is definable, via $\lim_{n\to\infty} (1+\tfrac1n)^n = \lim_{n\to\infty} \tfrac{(n+1)^n}{n^n}$;
more formally, $e$ is the real number $x$ satisfying the condition
\begin{equation*}
  \forall \varepsilon>0 \; \exists n\in\N\; \forall m\in\N \;\;
  \( m\ge n \imp -\varepsilon m^m < (m+1)^m-e m^m < \varepsilon m^m \).
\end{equation*}
This is not quite the definition mentioned in Introduction, but equivalent to it.

Similarly, for every convergent computable sequence of rational numbers, its limit is a definable number. In other
words, every \index{number!limit computable}\href{http://en.wikipedia.org/wiki/Computation in the
  limit\#Limit_computable_real_numbers}{limit computable} real number is definable.

Every computable real number is limit computable, therefore definable. In particular, the number $\pi$ is computable,
therefore definable.

Chaitin's constant is not computable, but still, limit computable (recall Introduction: it is the limit of a computable
increasing sequence of rational numbers), therefore definable. So, all the five constants discussed in Introduction
(taken from the book "Mathematical constants") are definable. Moreover, all the constants discussed in that book are
definable.

On the other hand, if we choose a number between 0 and 1 at random,\index{number!chosen at random} according to the
\href{http://en.wikipedia.org/wiki/Uniform distribution (continuous)}{uniform distribution}, we
\href{http://en.wikipedia.org/wiki/Almost surely}{almost surely} get an undefinable number, because the definable
numbers are a \href{http://en.wikipedia.org/wiki/Countable set}{countable set}.

Of course, such a randomly chosen undefinable number is not an explicit example of undefinable number. It may seem that
``explicit example of undefinable number'' is a patent nonsense, just as ``defined undefinable number''. But no, not quite
nonsense, see Section \ref{sect4}.

An infinite sequence $(x_1,x_2,\dots)=(x_n)_n$ of real numbers is nothing but the binary relation $\{(n,x)\mid n\in\N
\land x=x_n\} =$ $\{(1,x_1),(2,x_2),\dots\}$; if this binary relation is definable, we say that the sequence is
definable.\index{definable!sequence} If a sequence is definable, then all its members are definable numbers. However, a sequence of definable
numbers is generally not definable.

\begin{exercise}
 If a definable sequence converges, then its limit is a definable number. Prove it.
\end{exercise}

A function $f:\R \to \R$ is nothing but the binary relation ``$f(x)=y$'', that is, $A=\{(x,y)\mid f(x)=y\}$; if this
binary relation is definable, we say that the function is definable.\index{definable!function} An arbitrary binary relation $A$ is a function if
and only if for every $x$ there exists one and only one $y$ such that $(x,y)\in A$.

\begin{exercise}
 If $f$ is a definable function and $x$ is a definable number, then $f(x)$ is a definable number. Prove it.
\end{exercise}

However, a function that has definable values at all definable arguments is generally not definable.

\begin{exercise}
 If a definable function is differentiable, then its derivative is a definable function. Prove it.  \emph{Hint:} the
 derivative is the limit of\dots
\end{exercise}

\begin{exercise}
 If a definable function $f$ is continuous, then its
 \href{http://en.wikipedia.org/wiki/Antiderivative}{antiderivative} $F$ is definable if and only if $F(0)$ is a
 definable number. Prove it.  \emph{Hint:} $F(x)=F(0)+\lim_{n\to\infty} \frac x n \sum_{k=1}^n f(\frac k n x)$.
\end{exercise}

Similarly to the number $e$ we can treat the exponential function $x\mapsto e^x$. First, the relation $
\big\{(n,p,q,u)\mid n\in\N \land p\in\N \land q\in\N \land u=\(1+\frac p q \cdot \frac 1n\)^n\big\}$ is definable (since
$ \(1+\frac p q \cdot \frac 1n\)^n$ is a computable function of $n,p,q$). Second, the relation $\{(x,y)\mid y=e^x\}$ is
definable, since $e^x$ is the limit of $ \(1+\frac p q \cdot \frac 1n\)^n$ as $n$ tends to infinity and $\tfrac p q$
tends to $x$; more formally (but still not completely formally\dots), $y=e^x$ if and only if
\begin{multline*}
 \forall \varepsilon>0 \; \exists \delta>0 \; \forall n\in\N\, \forall p\in\Z\, \forall q\in\N\,
 \forall u \;\; \\
 \bigg( \Big(n\ge\frac1\delta\Big) \land \Big(-\delta<x-\frac p q<\delta\Big) \land \Big(u=\(1+\tfrac p q \cdot
 \tfrac 1n\)^n\Big) \imp \varepsilon < y-u < \varepsilon \bigg);
\end{multline*}
here $\Z$ is the set of integers (evidently definable).

The cosine function may be treated via complex numbers and \href{http://en.wikipedia.org/wiki/Euler's formula}{Euler's
  formula} $e^{ix}=\cos x + i \sin x$. First, the real part of the complex number $ \(1+i\frac p q \cdot \frac
1n\)^n$ is a computable function of $n,p,q$. Second, its limit as $n$ tends to infinity and $\tfrac p q$ tends to $x$
is equal to the real part $\cos x$ of the complex number $e^{ix}$.

Note that the \href{http://en.wikipedia.org/wiki/Exponential integral}{exponential integral} $\operatorname{Ei}(x)$
and the \href{http://en.wikipedia.org/wiki/Trigonometric integral\#Sine_integral}{sine integral} $\operatorname{Si}(x)$
are definable \href{http://en.wikipedia.org/wiki/Nonelementary integral}{nonelementary functions}.

Definable functions can be \href{http://en.wikipedia.org/wiki/Pathological (mathematics)}{pathological} and disrespect
dimension. In particular, there is a definable one-to-one correspondence between the (two-dimensional) square
$(0,1)\times(0,1)$ and a subset of the (one-dimensional) interval $(0,1)$, which will be used in Section \ref{sect6}.
Here is a way to this fact.

Given two numbers $x,y\in(0,1)$, we consider their decimal digits:\index{digit} $ x=(0.\al_1\al_2\dots)_{10}=\sum_{n=1}^\infty
10^{-n}\al_n$ where $\al_n\in\{0,1,2,3,4,5,6,7,8,9\}$ for each $n$, and the set $\{n:\al_n\ne9\}$ is infinite
(since we represent, say, $\tfrac12$ as $(0.5000\dots)_{10}$ rather than $(0.4999\dots)_{10}$); and similarly
$y=(0.\be_1\be_2\dots)_{10}$. We interweave their digits, getting a third number
$z=(0.\al_1\be_1\al_2\be_2\dots)_{10} \in (0,1)$. The ternary relation between such $x,y,z$ is a function
$W_2:(0,1)\times(0,1)\to(0,1)$.\index{zzW@$W_2$, injects 2-dim to 1-dim} Not all numbers of $(0,1)$ are of the form $W_2(x,y)$ (for example,
$\frac{21}{1100}=(0.01909090\dots)_{10}$ is not), which does not matter. It does matter that $x,y$ are uniquely
determined by $W_2(x,y)$, that is, $W_2(x_1,y_1)=W_2(x_2,y_2)$ implies $x_1=x_2 \land y_1=y_2$. In other words, $W_2$ is
an \href{http://en.wikipedia.org/wiki/Injective function}{injection}\index{injection} $(0,1)\times(0,1)\to(0,1)$.

Denoting by $D(n,x)$\index{zzD@$D(n,x)$, the $n$-th digit} the $n$-th decimal digit $\al_n$ of $x\in(0,1)$ we have $D(n,x)=\lfloor 10 \cdot
\operatorname{frac} (10^{n-1}x) \rfloor$; here $\lfloor a \rfloor$ is the \href{http://en.wikipedia.org/wiki/Floor and
  ceiling functions}{integer part} of $a$, and $\operatorname{frac}(a)=a-\lfloor a \rfloor$ is the
\href{http://en.wikipedia.org/wiki/Fractional part}{fractional part} of $a$.

\begin{exercise}
 The integer part function is definable. Prove it.  \emph{Hint:} $\{(x,\lfloor x \rfloor) \mid x>0\} = \{(x,n)\mid x>0
 \land n+1\in\N \land n\le x<n+1\}$.
\end{exercise}

\begin{exercise}
 The function $D:\N\times(0,1)\to\R$ is definable. (See Booij \cite[Lemma 3.4]{Bo}.) Prove it.  \emph{Hint:}
 $D(n,x)=d \equ \exists k\in\N \;\; \exists y\in\R \;\; \( k=10^{n-1} \land
 y=\operatorname{frac}(kx) \land d=\lfloor 10y\rfloor \)$.
\end{exercise}

\begin{exercise}
 The function $W_2:(0,1)^2\to(0,1)$ is definable. Prove it.  \emph{Hint:} $z=W_2(x,y) \equ \forall
 n\in\N \;\; \( D(2n,z)=D(n,y) \land D(2n-1,z)=D(n,x) \)$.
\end{exercise}

\begin{exercise}
 Generalize the previous exercise to $W_3:(0,1)^3\to(0,1)$. \emph{Hint:} consider $D(3n-2,z)$, $D(3n-1,z)$, $D(3n,z)$.
\end{exercise}

\section{Explicit example of undefinable number}
\label{sect4}

We construct such example in two steps. First, we enumerate all numbers definable in $(\R; +,\times,\N)$ (``first order''
is meant but omitted, as before); that is, we construct a sequence $(x_1,x_2,\dots)$ of real numbers that contains all
numbers definable in $(\R; +,\times,\N)$ (and only such numbers). Second, we construct a real number not contained in
this sequence.

The second step is well-known and simple, so let us do it now, for an arbitrary sequence $(x_1,x_2,\dots)$ of real
numbers. We construct a real number $x$ via its decimal digits, as $ x=\sum_{n=1}^\infty \frac{\al_n}{10^n}$, and we
choose each $\al_n$ to be different from the \nbdh{n}{th} digit (after the decimal point) of the absolute value
$|x_n|$ of $x_n$. To be specific, let us take $\al_n=3$ if $10k+7\le 10^n |x_n|<10k+8$ for some integer $k$, and
$\al_n=7$ otherwise. Then $x\ne x_n$ since the integral part of $10^n|x|$, being of the form $10\ell+\al_n$ for integer
$\ell$, is different from the integral part of $10^n|x_n|$, the latter being of the form $10k+\be_n$ for integer $k$,
and $\al_n\ne\be_n$ (either $\be_n=7,\al_n=3$ or $\be_n\ne7,\al_n=7$). This is an instance of
\href{http://en.wikipedia.org/wiki/Cantor's diagonal argument}{Cantor's diagonal argument}.\index{Cantor's diagonal argument}

Now we start constructing a sequence $(x_1,x_2,\dots)$ of real numbers that contains all numbers definable in $(\mathbb
R; +,\times,\N)$ (and only such numbers). These numbers being elements of single-element subsets of $\R$ definable in
$(\R; +,\times,\N)$, and these subsets being unary relations, we enumerate all relations (unary, binary, \dots) definable
in $(\R; +,\times,\N)$. These are obtained from the three given relations (addition, multiplication, ``naturality'') via
the 5 operations (complement, union, permutation, set multiplication, projection) applied repeatedly. We may save on
permutations by restricting ourselves to \href{http://en.wikipedia.org/wiki/Cyclic_permutation\#Transpositions}{adjacent
  transpositions},\index{transposition, adjacent} that is, permutations that swap two adjacent numbers $k,k+1$ and leave intact other numbers of
$\{1,\dots,n\}$; this is sufficient, since every permutation is a product of some adjacent transpositions. We start with
the three given relations
\begin{align*}
A_1&=\{(x,y,z)\mid x+y=z\}, &&\text{``addition''} \\
A_2&=\{(x,y,z)\mid xy=z\}, &&\text{``multiplication''} \\
A_3&=\N, &&\text{``naturality''}
\end{align*}
and apply to them the five operations (whenever possible). The first operation ``complement'' gives
\begin{align*}
A_4&=\{(x,y,z)\mid x+y\ne z\}, \\
A_5&=\{(x,y,z)\mid xy\ne z\}, \\
A_6&=\{x\mid x\notin\N\}.
\end{align*}
The ``union'' operation gives
\begin{equation*}
A_7=\{(x,y,z)\mid (x+y=z) \lor (xy=z)\}.
\end{equation*}
The ``permutation'' operation (reduced to adjacent transpositions), applied to the ternary relation $A_1$, gives two
relations $A_8,A_9$; namely, $A_8=\{(x,y,z)\mid y+x=z\}$ (equal to $A_1$ due to commutativity, but we do not bother) and
$A_9=\{(x,y,z)\mid x+z=y\}$; we apply the same to $A_2$ getting $A_{10},A_{11}$. Further, ``set multiplication'' gives the
4-ary relation
\begin{equation*}
A_{12}=\{(x,y,z,w)\mid x+y=z\},
\end{equation*}
similarly $A_{13}$, and $A_{14}=\{(x,y)\mid x\in\N\}$. The most remarkable ``projection'' operation gives
\begin{equation*}
  A_{15}=\{(x,y)\mid \exists z \; (x+y=z)\}
\end{equation*}
(in fact, $A_{15}=\R^2$), and similarly $A_{16}$.

The first 16 relations $A_1,\dots,A_{16}$ are thus constructed. On the next iteration we apply the 5 operations to these
16 relations (whenever possible; though, some are superfluous) and get a longer finite list. And so on, endlessly. A bit
cumbersome, but really, a routine exercise in programming, isn't it? Well, it is, provided however that the
``programming language'' stipulates the data type ``relation over $\R$'' and the relevant operations on relations. By
the way, equality test for relations is not needed (unless we want to skip repetitions); but test of existence and
uniqueness (for unary relations), and extraction of the unique element, are needed for the next step.

Now we are in position to construct $x_n$; for each $n$ we check, whether the relation $A_n$ is of the form $\{u\}$ for
$u\in\R$ or not; if it is, we take $x_n=u$, otherwise $x_n=0$. (Note that $x_n=0$ whenever the relation $A_n$ is not
unary.)

Applying the diagonal argument (above) to this sequence $(x_1,x_2,\dots)$ we construct a real number $x$ not contained
in the sequence, therefore, not definable in $(\R; +,\times,\N)$.

This number $x$ is defined, but not in $(\R; +,\times,\N)$. Why not? Because the definition of $x$ involves a sequence
of relations in $\R$. Sequences of numbers are used in Section \ref{sect3}, but sequences of relations are something
new, beyond the first order. (See Wikipedia: \href{http://en.wikipedia.org/wiki/First-order logic}{First-order logic},
\href{http://en.wikipedia.org/wiki/Second-order logic}{Second-order logic}.)

Is there a better approach? Could we define in $(\R; +,\times,\N)$ the same sequence $(x_1,x_2,\dots)$, or maybe another
sequence containing all computable numbers, by a clever trick? No, this is impossible. For every sequence definable in
$(\R; +,\times,\N)$ the diagonal argument gives a number definable in $(\mathbb R; +,\times,\N)$ and not contained in
the given sequence.

\begin{quote}
$\bullet$ \emph{there is no definable enumeration of definable reals} (Poincar\'e 1909), see
  \href{https://plato.stanford.edu/entries/paradoxes-contemporary-logic/}{Stanford Encyclopedia of Philosophy: Paradoxes
    and Contemporary Logic}.
\end{quote}

\section{Second order}
\label{sect5}

We introduce second-order definability in $(\R; +, \times)$. The set $\N$ of natural numbers is second-order definable
in $(\R; +, \times)$, as we'll see soon. In contrast to the first order definability, usual definitions of mathematical
constants will apply without recourse to computability and Diophantine sets.

A \href{http://en.wikipedia.org/wiki/Second-order predicate}{second-order predicate} is a predicate that takes a
first-order predicate as an argument. Likewise, a second-order relation\index{relation!second-order} is a relation between relations. For example,
the binary relation $f'=g$ between a function $f$ and its derivative $g$ may be thought of as a relation between two
binary relations: first, the relation $f(x)=y$ between real numbers $x,y$, and second, the relation $g(t)=v$ between
real numbers $t,v$. First order definability of a real number involves definable first-order relations (between real
numbers). Second order definability of a real number involves definable second-order relations (between first-order
relations). Here is a possible formalization of this idea.

We introduce the set\index{zzS@$\mathbf S$, set of all tuples and relations}\index{zzP@$\oP(\R^n)$, set of relations}
\begin{equation*}
 \! \mathbf S = \( \R \cup \R^2 \cup \R^3 \cup \dots \) \cup \( \oP(\R)
 \cup \oP(\R^2) \cup \oP(\R^3) \cup \dots \) = \bigg( \bigcup_{n=1}^\infty
 \R^n \bigg) \cup \bigg( \bigcup_{n=1}^\infty \oP(\R^n) \bigg)
\end{equation*}
that contains, on one hand, all tuples\index{tuple} (finite sequences) $(x_1,\dots,x_n)\in\R^n$ of real numbers (for all $n$; here we
do not distinguish 1-tuples from real numbers), and on the other hand, all $n$-ary relations $A\subset\mathbb R^n$ on
$\R$ (for all $n$). Here $\oP(\R)$\index{zzP@$\oP(\R)$, power set} is the set of all subsets (that is, the \href{http://en.wikipedia.org/wiki/Power
  set}{power set}) of the real line $\R$, in other words, of unary relations on $\R$; $\oP(\R^2)$ is the set of all
subsets (that is, the \href{http://en.wikipedia.org/wiki/Power set}{power set}) of the Cartesian plane $\R^2$, in other
words, of binary relations on $\R$; and so on. On this set $\mathbf S$ we introduce two relations:
\begin{compactitem}
\item membership,\index{relation!membership} the binary relation $\cup_{n=1}^\infty \{ \((x_1,\dots,x_n),A\)\mid (x_1,\dots,x_n)\in A \}
\subset \mathbf S^2 $; it says that the given $n$-tuple belongs to the given $n$-ary relation;
\item ``appendment'',\index{relation!appendment} the ternary relation $\cup_{n=1}^\infty \{ \( x, (x_1,\dots,x_n), (x_1,\dots,x_n,x) \mid
x_1,\dots,x_n,x \in \R \} \subset \mathbf S^3 $; it says that the latter tuple results from the former tuple by
appending the given real number.
\end{compactitem}

The two ternary relations on $\R$, addition and multiplication, may be thought of as ternary relations on
$\mathbf S$ (since $\R \subset \mathbf S$):
\begin{compactitem}
\item addition: $\{(x,y,z)\mid x,y,z\in\R, x+y=z\} \subset \mathbf S^3$;
\item multiplication: $\{(x,y,z)\mid x,y,z\in\R, xy=z\} \subset \mathbf S^3$.
\end{compactitem}

We endow $\mathbf S$ with the \Dstructure\ generated by the four relations
(membership, appendment, addition, multiplication). All relations on $\mathbf S$ that belong to this
\Dstructure\ will be called second-order definable.\index{second-order definable} In the rest of this section,
``definable'' means ``second-order definable'', unless stated otherwise.

\begin{exercise}
 The set $\cup_{n=1}^\infty \R^n$ of all tuples and the set $\cup_{n=1}^\infty \oP(\R^n)$ of
 all relations are definable subsets of $\mathbf S$. Prove it.  \emph{Hint:} first, the set of all tuples is $\{ s\in
 \mathbf S \mid \exists r,t \in \mathbf S \; (r,s,t)\in A \}$ where $A$ is the appendment relation; second, take the
 complement.
\end{exercise}

\begin{exercise}
 The set $\R$ of all real numbers is a definable subset of $\mathbf S$. Prove it.  \emph{Hint:} $\R = \{
 r\in \mathbf S \mid \exists s,t\in \mathbf S \;(r,s,t)\in A \}$ where $A$ is the appendment relation.
\end{exercise}

\begin{exercise}
 Each $\R^n$ is a definable subset of $\mathbf S$. Prove it.  \emph{Hint:} $\R^{n+1} = \{ t\in \mathbf S
 \mid \exists r \in \R \; \exists s \in \R^n \; (r,s,t)\in A \}$ where $A$ is the appendment relation.
\end{exercise}

\begin{exercise}
 The set $\oP(\R)$ of all sets of real numbers (that is, unary relations) is a definable subset of $\mathbf S$. Prove
 it.  \emph{Hint:} for $A\in\cup_{n=1}^\infty \oP(\R^n)$ we have $A\in \oP(\R) \equ A\subset\R \equ ( \forall a\in A
 \;\> a\in\R) \equ \neg ( \exists a\in A \;\> a\notin\R)$; apply the projection to $\{(A,a) \mid A \in \cup_{n=1}^\infty
 \oP(\R^n) \land a\in A \land a \notin \R \}$.
\end{exercise}

\begin{exercise}
 For each $n$ the set $\oP(\R^n)$ of all $n$-ary relations is a definable subset of $\mathbf S$. Prove it.
 \emph{Hint:} similar to the previous exercise.
\end{exercise}

\begin{exercise}
 If $B\subset \oP(\R)$ is a definable set of subsets of $\R$, then the union $\cup_{A\in B} A$ of all these subsets is a
 definable set (of real numbers). Prove it.  \emph{Hint:} for $x\in\R$ we have $x\in\cup_{A\in B} A \equ \exists A \;
 (A\in B \land x\in A)$; take the projection of $\{(x,A) \mid A\in B \land x\in A \} = (\R \times B) \cap \{(x,A)\mid
 x\in A\}$.
\end{exercise}

\begin{exercise}
 Do the same for the intersection $\cap_{A\in B} A$. \emph{Hint:} $\cap_{A\in B} A = \R \setminus \cup_{A\in B} (\R
 \setminus A)$; consider $\{(x,A) \mid A\in B \land x\in \R \setminus A \} = (\R \times B) \cap \{(x,A)\mid x\notin
 A\}$. (But what if $B$ is empty?)
\end{exercise}

\begin{exercise}
 Generalize the two exercises above to $B\subset \oP(\R^2)$, $B\subset \oP(\R^3)$ and so on. \emph{Hint:} now $x$ is a
 tuple.
\end{exercise}

In particular, taking a single-element set $B=\{A\}$ we see that definability of $\{A\}$ implies definability of
$A$. The converse holds as well (see below).

\enlargethispage{10pt}

\begin{exercise}
If a set $A\subset\R$ (of real numbers) is definable, then the set $\oP(A)$ (of all subsets of $A$) is definable. Prove
it.  \emph{Hint:} for $A_1\in\oP(\R)$ we have $A_1\in\oP(A) \equ A_1 \subset A \equ (\forall x\in A_1\;\> x\in A)\equ \neg
( \exists x\in A_1\;\> x\notin A)$; consider $\{(A_1,x) \mid A_1\in\oP(\R) \land x\in\R \land x\in A_1 \land x\notin A\} =
\( \oP(\R) \times (\R\setminus A) \) \cap \{(A_1,x)\mid x\in A_1\}$.
\end{exercise}%

\begin{exercise}
 Do the same for the set $\{A_1\in \oP(\R) \mid A_1\supset A\}$ (of all supersets of $A$).  \emph{Hint:} similarly to
 the previous exercise, consider $\( P(\R)\times A\) \cap \{(A_1,x)\mid x\notin A_1\}$.
\end{exercise}

\begin{exercise}
 Generalize the two exercises above to $A\subset\R^2$, $A\subset\R^3$ and so on.
\end{exercise}

\vspace{-3pt}

\emph{Remark.} These 11 exercises (above) do not use addition and multiplication, nor any properties of real
numbers. They generalize readily to a more general situation. One may start with an arbitrary set $R$ (rather than the
real line $\R$), consider the set $S$ constructed from $R$ as above (all tuples and all relations), endow $S$ by a
\Dstructure\ such that the two relations on $S$, membership and appendment, are definable, and
generalize the 11 exercises to this case.

Taking the intersection of the set of subsets and the set of supersets we see that definability of $A$ implies
definability of the single-element set (called \href{http://en.wikipedia.org/wiki/Singleton (mathematics)}{singleton})
$B=\{A\}$. So, $A$ is definable if and only if $\{A\}$ is definable. And still (by convention, as before) a real number
$x$ is definable if and only if $\{x\}$ is definable.

Does it mean that, for example, numbers $0$, $1$, $\sqrt2$ are definable (as well as every rational number and every
algebraic number)? We know that they are first-order definable in $(\R; +,\times)$; does it follow that they are
(second-order) definable in $\mathbf S?$

The answer is affirmative, but needs a proof. Here we face another general question. Let $S$ be a set and $R\subset S$
its subset. Every \ary{n} relation on $R$ is also a \ary{n} relation on
$S$ (since $R\subset S \imp R^n\subset S^n \imp \oP(R^n) \subset \oP(S^n)$). Thus, given some relations on $R$, we get
two \Dstructure s; first, the \Dstructure\ on $R$ generated by the given
relations, and second, the \Dstructure\ on $S$ generated by the same relations.

\enlargethispage{3pt}

\textbf{Lemma.}  Assume that $R$ is a definable subset of $S$ (according to the second \Dstructure). Then every relation
on $R$ definable according to the first \Dstructure\ is also definable according to the second \Dstructure.\footnote{%
  \emph{Proof.} Denote the first \Dstructure\ by $D_R$ and the second by $D_S$.  We know that $R\in D_S$. It follows
  (via set multiplication) that $R\times S \in D_S$, $R\times S\times S = R\times S^2 \in D_S$, and so on; by induction,
  $R\times S^n \in D_S$ for all $n$. Thus (via permutation), $S^n\times R \in D_S$.
  
  In order to prove that $D_R \subset D_S$ we compare the five operations on relations (complement, union, permutation,
  set multiplication, projection) over $R$ (call them \nbdh{R}{operations}) and over $S$
  (\nbdh{S}{operations}). We have to check that each \nbdh{R}{operation} applied to
  relations on $R$ that belong to $D_S$ gives again a relation (on $R$) that belongs to $D_S$.

  For the union we have nothing to check, since the \nbdh{R}{union} of two relations is equal to their
  \nbdh{S}{union}. Similarly, we have nothing to check for permutation and projection. Only set
  multiplication and complement need some attention.

  Set multiplication. The \nbdh{R}{multiplication} applied to $A\in\oP(R^n)\cap D_S$ gives
  $A\times R$. We have $ A\times R = (A\times S) \cap (S^n\times R) \in D_S$ since $A\times S \in D_S$ and $S^n\times R
  \in D_S$.

  It follows (by induction) that $R^n\in D_S$ for all $n$.

  Complement. The \nbdh{R}{complement} applied to $A\in\oP(R^n)\cap D_S$ gives
  $R^n\setminus A$. We note that the \nbdh{S}{complement} $S^n\setminus A$ belongs to $D_S$
  (since $A\in D_S$), thus $R^n\setminus A = R^n \cap (S^n\setminus A) \in D_S$ (since $R^n\in D_S$).}

Now we are in position to prove definability of the set $\N$ of natural numbers. It is sufficient to prove definability
of the set $B\subset\oP(\R)$ of all sets $A\subset\R$ satisfying the two conditions $1\in A$ and $\forall x\in A \;
(x+1\in A)$ (since the intersection of all these $A$ is $\N$). The complement $\oP(\R)\setminus B=\{A\in\oP(\R) \mid
\exists x\in\R \; (x\in A\land x+1\notin A)\}$ is the projection of the intersection of two sets, $\{(A,x)\in\oP(\mathbb
R)\times\R \mid x\in A\}$ and $\{(A,x)\in\oP(\R)\times\R \mid x+1\notin A\}$. The former results from the (permuted)
membership relation; the latter is the projection of the projection of $\{(A,x,y,z)\in \oP(\R) \times \R^3 \mid
y\in\{1\} \land x+y=z \land z\notin A\}$, this set being the intersection of three sets: first, $\oP(\R) \times \R\times
\{1\} \times \R$; second, $\oP(\R)$ times the addition relation; third, $\oP(\R) \times \mathbb R\times\R \times
(\R\setminus A)$. It follows that $B$ is definable, whence $ \N = \cap_{A\in B} A$ is definable.

This is instructive. \emph{In order to formalize a definition of a set via its defining property, we have to deal with
  sets of sets, and more generally, relations between sets.}

Using again the lemma above we see that all real numbers first-order definable in $(\R; +,\times,\N)$ are second-order
definable. Section \ref{sect3} gives many examples, including the five numbers $ \sqrt {2}, \varphi, e, \pi, \Om, $
discussed in Introduction. But second-order proofs of their definability are much more easy and natural.

The binary relation ``$x\in\N \land y=x!$'' is the sequence $(n!)_{n\in\N}$ of
\href{http://en.wikipedia.org/wiki/Factorial}{factorials}, that is, the set
$\{(1,1),(2,2),(3,6),(4,24),(5,120),\dots\}$. It is definable, similarly to $\N$, since it is the \emph{least} subset
$A$ of $\R^2$ such that $(1,1)\in A$ and $(x,y)\in A \imp \(x+1,(x+1)y\)\in A$.  Alternatively, it is definable since it
is the \emph{only} subset $A$ of $\R^2$ with the following three properties:
\begin{gather*}
 \forall (x,y)\in A \;\; x\in\N, \\
 \forall x\in\N \;\; \exists! y\in\R \;\; (x,y)\in A, \\
 \forall (x,y)\in A \;\; \(x+1,(x+1)y\)\in A.
\end{gather*}
That is, the factorial is the only function $\N \to \R$ satisfying the
\href{http://en.wikipedia.org/wiki/Recurrence relation}{recurrence relation} $(n+1)!=(n+1)n!$ and the initial
condition $1!=1$.

\begin{exercise}
 Partial sums of the series $ \sum_{n=0}^\infty \frac1{n!}$ are a definable sequence. Prove it.
\end{exercise}

\begin{exercise}\label{5.13}
 The number $e$ is definable. Deduce it from the previous exercise.
\end{exercise}

In the first-order framework it is possible to treat many functions (for instance, the exponential function $x\mapsto
e^x$, the sine and cosine functions $\sin, \cos$, the exponential integral $\operatorname{Ei}$ and the sine integral
$\operatorname{Si}$) and many relations between functions (for instance, derivative and antiderivative); arguments and
values of these functions are arbitrary real numbers (not necessarily definable), but the functions are definable. Such
notions as arbitrary functions (not necessarily definable), continuous functions (and their antiderivatives),
differentiable functions (and their derivatives) need the second-order framework.

As was noted there, a function $f:\R \to \R$ is nothing but the binary relation ``$f(x)=y$'', that is,
$A=\{(x,y)\mid f(x)=y\}$. An arbitrary binary relation $A$ is such a function if and only if for every $x$ there exists
one and only one $y$ such that $(x,y)\in A$ (existence and uniqueness). For functions defined on arbitrary subsets of
the real line the condition is weaker: for every $x$ there exists at most one $y$ such that $(x,y)\in A$ (uniqueness).

\begin{exercise}
 (a) All $A\in\oP(\R^2)$ satisfying the uniqueness condition are a definable subset of $\oP(\R^2)$; (b) the same holds
  for the existence and uniqueness condition. Prove it.
\end{exercise}

\begin{exercise}
 All continuous functions $\R \to \R$ are a definable subset of $\oP(\R^2)$. Prove it.
\end{exercise}

\begin{exercise}
 All differentiable functions $\R \to \R$ are a definable subset of $\oP(\R^2)$. Prove it.
\end{exercise}

\begin{exercise}
 The binary relation ``$f'=g$'' is definable. That is, the set of all pairs $(f,g)$ of functions $\R \to \R$ such that
 $\forall x\in\R \; \(f'(x)=g(x)\)$ is a definable subset of $\oP(\R^2)\times \oP(\R^2)$ (in other words, definable
 binary relation on $\oP(\R^2)$). Prove it.
\end{exercise}

Antiderivative can now be treated in full generality. In contrast, in the first-order framework it was treated via
Riemann integral $ F(x)=F(0)+\lim_{n\to\infty} \frac x n \sum_{k=1}^n f(\frac k n x)$ for continuous definable $f$
only. In particular, now the exponential function $x\mapsto e^x$ may be treated via $f(e^x)-f(1)=x$ where
$f'(x)=\tfrac1x$ for $x>0$; accordingly, the constant $e$ may be treated via $f(e)-f(1)=1$. Alternatively, the
exponential function may be treated via the differential equation $f'=f$ (and initial condition $f(0)=1$). Trigonometric
functions $\sin, \cos$ may be treated via the differential equation $f''=-f$; accordingly, the constant $\pi$ may be
treated as the least positive number such that $(f''=-f) \imp \(f(\pi)=-f(0)\)$. Or, alternatively, as $ \pi = 4\int_0^1
\sqrt{1-x^2} \, dx$ (via antiderivative).

This is instructive. \emph{In the second-order framework we may define functions (and infinite sequences) via their
  properties,} irrespective of computability, Diophantine equations and other tricks of the first-order framework.

Nice; but what about second-order definable real numbers? Are they all first-order definable, or not? Even if obtained
from complicated differential equations, they are computable, therefore, first-order
definable in $(\R; +,\times,\N)$. Probably, our only chance to find a second-order definable but first-order undefinable
number is, to prove that the explicit example of (first-order) undefinable number, given in Section~\ref{sect4}, is
second-order definable; and our only chance to prove this conjecture is, to formalize that section within the
second-order framework.

\section{First-order undefinable but second-order definable}
\label{sect6}

Recall the infinite sequence of relations $(A_k)_{k=1}^\infty$ treated in Section \ref{sect4}. Is it second-order
definable? Each $A_k$ belongs to the set $\mathbf S$ (from Section~\ref{sect5}); their infinite sequence is a binary
relation between $k$ and $A_k$ (namely, the set of pairs $\{(1,A_1),(2,A_2),\dots\}$), thus, a special case of a binary
relation on $\mathbf S$; the question is, whether this relation is definable, or not. Like the sequence of factorials,
it is defined by recursion. But factorials, being numbers, are first-order objects, which is why their sequence is
second-order definable via its properties. In contrast, relations $A_k$ are second-order objects! Does it mean that
third order is needed for defining their sequence by recursion?

True, the sequence of factorials is first-order definable (over $(\R; +,\times,\N)$) due to its computability, via
Matiyasevich's theorem. Could something like that be invented for second-order objects? Probably not.

Yet, these obstacles are surmountable. The sequence of relations may be replaced with a single relation by a kind of
\href{http://en.wikipedia.org/wiki/Currying}{currying} (or rather, uncurrying); the
\index{disjoint union}\href{http://en.wikipedia.org/wiki/Disjoint union}{disjoint union} $\{1\}\times A_1 \cup \{2\}\times A_2 \cup \dots$
may be used instead of the set of pairs $\{(1,A_1),(2,A_2),\dots\}$. Further, relations $A_k$ of different
\href{http://en.wikipedia.org/wiki/Arity}{arities} may be replaced with unary relations (subsets of the real line),
since two real numbers may be encoded into a single real number via an appropriate definable
\href{http://en.wikipedia.org/wiki/Injective function}{injection} $\R^2 \to \R$, and the same applies to
three and more numbers (moreover, to infinitely many numbers, see Booij \cite[Sect. 3.2]{Bo}). In addition, tuples
$(x_1,\dots,x_n)$ may be replaced (whenever needed) by finite sequences $(x_k)_{k=1}^n = \{(1,x_1),\dots,(n,x_n)\}$,
which provides a richer assortment of definable relations.

The distinction between tuples and finite sequences\index{tuple!versus finite sequence} is a technical subtlety that may be ignored in many contexts, but
sometimes requires attention. It is tempting to say that an ordered pair $(a,b)$, a 2-tuple (that is, tuple of length
2), and a 2-sequence (that is, finite sequence of length 2) are just all the same. However, the 2-sequence is, by
definition, a function on $\{1,2\}$, thus, the set of two ordered pairs $\{(1,a),(2,b)\}$. Surely we cannot define an
ordered pair to be a set of two other ordered pairs! If sequences are defined via functions, and functions are defined
via pairs, then pairs must be defined before sequences, and cannot be the same as 2-sequences. See Wikipedia:
\href{http://en.wikipedia.org/wiki/Sequence\#Formal_definition_and_basic_properties}{sequence (formal definition)},
\href{http://en.wikipedia.org/wiki/Tuple\#Definitions}{tuples (as nested ordered pairs)}, and
\href{http://en.wikipedia.org/wiki/Ordered_pair\#Defining_the_ordered_pair_using_set_theory}{ordered pair: Kuratowski's
  definition}. For convenience we'll denote a finite sequence $(x_k)_{k=1}^n = \{(1,x_1),\dots,(n,x_n)\}$ by
$[x_1,\dots,x_n]$; it is similar to, but different from, the tuple $(x_1,\dots,x_n)$.

We'll construct again, this time in the second-order framework, the sequence $(x_1,x_2,\dots)$ of real numbers that
contains all numbers first-order definable in $(\R; +,\times,\N)$, exactly the same sequence as in
Section \ref{sect4}. To this end we'll construct first the disjoint union $\{1\}\times B_1
\cup \{2\}\times B_2 \cup \dots$ of unary relations $B_k$ on $\R$ similar to, but different from, relations
$A_1,A_2,\dots$ (unary, binary, ...) constructed there (that exhaust all relations first-order definable in $(\mathbb
R; +,\times,\N)$).

Before the unary relations $B_k$ we construct 4-tuples $b_k$ of integers (call them "instructions") imitating a program
for a machine that computes $B_k$. Similarly to a \href{http://en.wikipedia.org/wiki/Machine code}{machine language
  instruction}, each $b_k$ contains an \href{http://en.wikipedia.org/wiki/Opcode}{operation code}, address of the
first \href{http://en.wikipedia.org/wiki/Operand\#Computer_science}{operand}, a parameter or address of the second
operand (if applicable, otherwise 0), and in addition, the arity of $A_k$.

Recall Section \ref{sect4}. Three relations $A_1,A_2,A_3$ of arities $3,3,1$ are given, and lead to the next $13$
relations $A_4,\dots,A_{16}$. In particular, $A_4$ is the complement of $A_1$. Accordingly, we let $b_4=(1,1,0,3)$;
here, operation code $1$ means ``complement...'', operand address $1$ means ``...of $A_1$'', the third number $0$ is
dummy, and the last number $3$ means that the relation $A_4$ is ternary. Similarly, $b_5=(1,2,0,3)$ and $b_6=(1,3,0,1)$.

Further, $A_7$  being the union  of $A_1$ and  $A_2$, we let $b_7=(2,1,2,3)$;  operation code $2$  means ``union\dots'',
first operand address $1$ means ``\dots of $A_1$'', second operand address $2$ means ``\dots and $A_2$'', and again, $3$
is the arity of $A_7$.

Further, $A_8$ being a permutation of $A_1$, we let $b_8=(3,1,1,3)$; operation code $3$ means ``permutation\dots'',
operand address 1 means ``\dots of $A_1$'', the parameter $1$ means ``swap $1$ and $2$'', and $3$ is the arity of $A_8$.
Similarly, $b_9=(3,1,2,3)$ (in $A_1$ swap $2$ and $3$), $b_{10}=(3,2,1,3)$ (in $A_2$ swap $1$ and $2$),
$b_{11}=(3,2,2,3)$ (in $A_2$ swap $2$ and $3$).

Further, $A_{12}$ being $A_1\times\R$, we let $b_{12}=(4,1,0,4)$; operation code $4$ means ``set multiplication'', $1$
refers to the operand $A_1$, and $4$ is the arity of $A_{12}$. Similarly, $b_{13}=(4,2,0,4)$ and $b_{14}=(4,3,0,2)$.

Further, $A_{15}$ being the projection of $A_1$, we let $b_{15}=(5,1,0,2)$; operation code 5 means "projection\dots",
$1$ means "\dots of $A_1$", and $2$ means "\dots is binary". Similarly, $b_{16}=(5,2,0,2)$.

This way, the finite sequence $[3,3,1]$ of natural numbers (interpreted as arities) leads to the finite sequence
$[b_4,\dots,b_{16}]$ of 4-tuples (interpreted as instructions). Similarly, every finite sequence of natural numbers
leads to the corresponding finite sequence of 4-tuples. The relation between these two finite sequences is definable;
the proof is rather cumbersome, like a routine exercise in programming, but doable. Having this relation, we define an
infinite sequence of 4-tuples $b_k$ (interpreted as the infinite ``program'') together with an infinite sequence
$(k_n)_{n=1}^\infty$ by the following defining properties:
\begin{compactitem}
  \item $k_1=3; \;\;\; \forall n\in\N \;\; k_n<k_{n+1}$;
  \item $b_1=b_2=(0,0,0,3); \; b_3=(0,0,0,1)$;
  \item for every $n=1,2,\dots$ the finite sequence $[b_{k_n+1},\dots,b_{k_{n+1}}]$ of 4-tuples corresponds (according
    to the definable relation treated above) to the finite sequence of the natural numbers that are the last (fourth)
    elements of the 4-tuples $b_1,\dots,b_{k_n}$.
\end{compactitem}
In particular, $k_1=3$, $k_2=16$; the third property for $n=1$ states that $[b_4,\dots,b_{16}]$ corresponds to
$[3,3,1]$. And for $n=2$ it states that $[b_{17},\dots,b_{k_3}]$ corresponds to $[3,3,1,3,3,1,3,3,3,3,3,4,4,2,2,2]$. And
so on.

The infinite program is ready. It could compute all relations $A_k$ if executed by a machine able to process relations
of all arities. Is such machine available in our framework? The disjoint union $\{1\}\times A_1 \cup \{2\}\times A_2
\cup \dots$ could be used instead of the set of pairs $\{(1,A_1),(2,A_2),\dots\}$, but is not contained in (say)
$\R^{100}$. True, in practice 100-ary relations do not occur in definitions; but we investigate definability in
principle (rather than in practice). We encode all relation into unary relations as follows.

We recall the definable injective functions $W_2:(0,1)^2\to(0,1)$ and $W_3:(0,1)^3\to(0,1)$ treated in the end of
Section \ref{sect3}. The same works for any $(0,1)^m$. But we need to serve all dimensions $m$ by a single definable
function. To this end we turn from sets $\R^m$ of $m$-tuples $(x_1,\dots,x_m)$ to sets, denote them $\mathbb
R^{[m]}$, of $m$-sequences $[x_1,\dots,x_m]$.

\begin{exercise}\label{6.1}
 The set of all finite sequences of real numbers, $\{ [x_1,\dots,x_m] \mid m\in\N,\,x_1,\dots,x_m\in\R\} =
 \cup_{m=1}^\infty \R^{[m]} \subset \oP(\R^2)$, is definable, and the binary relation "length", $\{ ([x_1,\dots,x_m],m) \mid
 m\in\N,\,x_1,\dots,x_m\in\R\} = \cup_{m=1}^\infty \( \R^{[m]} \times \{m\} \) \subset \oP(\R^2)\times\R$, on $\mathbf
 S$ is a definable function on that set. Prove it.  \emph{Hint:} start with the binary relation.
\end{exercise}

\begin{exercise}\label{6.2}
 The function $E:([x_1,\dots,x_m],k)\mapsto x_k$\index{zzE@$E$, evaluation} ("evaluation") is a definable real-valued function on the set
 $\cup_{m=1}^\infty \( \R^{[m]} \times \{1,\dots,m\} \)$. Prove it.  \emph{Hint:} for $s=[x_1,\dots,x_m]$, $k\in\N$ and
 $x\in\R$ we have $E(s,k)=x \equ (k,x)\in s \equ \exists p \( p\in s \land (k,x)=p \)$; use the two given relations on
 $\mathbf S$, membership and appendment (consider $n=1$ in the definition of appendment).
\end{exercise}

We choose a definable bijection $h:\R\to(0,1)$, for example, $ h(x)=\frac12(1+\frac{x}{1+|x|})$, and define a function
$W:\cup_{m=1}^\infty \R^{[m]} \to \R$\index{zzW@$W$, injects $n$-dim to 1-dim} by $W([x_1,\dots,x_m])=W_m\(h(x_1),\dots,h(x_m)\)$ for all $m\in\N$ and
$x_1,\dots,x_m\in\R$.

\begin{sloppypar}
\begin{exercise}
 The function $W$ is definable. Prove it.  \emph{Hint:} for all $m\in\N$, $x_1,\dots,x_m\in\R$ and $z\in\R$ we have
 $W([x_1,\dots,x_m])=z \equ W_m\(h(x_1),\dots,h(x_m)\)=z \equ \forall k\in\{1,\dots,m\} \;\; \forall n\in\N \;\;
 D(m(n-1)+k,z) = D(n,h(x_k))$; that is, for all $m\in\N$, $s\in\R^{[m]}$ and $z\in\R$ holds $W(s)=z \equ \forall k\in\N
 \;\; \( k\le m \imp \forall n\in\N \;\; D(m(n-1)+k,z)=D(n,h(E(s,k)) \)$.
\end{exercise}
\end{sloppypar}

At last, we are in position to ``execute the infinite program'' $(b_k)_{k=4}^\infty$, that is, to prove (second-order)
definability of the set $B = \{1\}\times B_1 \cup \{2\}\times B_2 \cup \dots \subset \R^2$, the disjoint union of unary
relations $B_k$ on $\R$ that encode (according to $W$) the relations $A_1,A_2,\dots$ (that exhaust all relations
first-order definable in $(\R; +,\times,\N)$).

We extract $B_1=\{c\in\R \mid (1,c)\in B\}$, decode the ternary relation $\{[x,y,z]\mid W([x,y,z])\in B_1\} =
\{s\in\R^{[3]} \mid W(s)\in B_1\}$ and require it to be (like $A_1$) the addition relation $\{[x,y,z] \mid x+y=z\} =
\{s\in \R^{[3]} \mid E(s,1)+E(s,2)=E(s,3)\}$. That is, we require
\begin{equation*}
 \forall s\in\R^{[3]} \;\; \( (1,W(s))\in B \equ E(s,1)+E(s,2)=E(s,3) \).
\end{equation*}
This condition fails to uniquely determine the set $B_1$, since the image of $\R^{[3]}$ under $W$ is not the whole $\R$
(not even the whole $(0,1)$). We prevent irrelevant points by requiring in addition that $\forall x\in\R \;\; \(
(1,x)\in B \imp \exists s\in \R^{[3]} \;\; W(s)=x \)$. We do not repeat such reservation below.

Similarly, $B_2$ must encode the multiplication relation, and $B_3$ must encode the set of natural numbers:
\begin{gather*}
 \forall s\in\R^{[3]} \;\; \( (2,W(s))\in B \equ E(s,1) E(s,2)=E(s,3) \); \\
 \forall s\in\R^{[1]} \;\; \( (3,W(s))\in B \equ E(s,1) \in \N \).
\end{gather*}
These first three requirements (above) are special. Other requirements should be formulated in general, like this: for
every $k\ge4$, if the first element of $b_k$ (the operation code) equals 1, then (...), otherwise (...). But let us
consider several examples before the general case.

According to the instruction $b_4$, the set $B_4$ must encode the complement of the set encoded by $B_1$:
\begin{equation*}
 \forall s\in\R^{[3]} \;\; \( (4,W(s))\in B \equ (1,W(s))\notin B \);
\end{equation*}
similarly,
\begin{align*}
 \forall s\in\R^{[3]} \;\; &\( (5,W(s))\in B \equ (2,W(s))\notin B \); \\
 \forall s\in\R^{[1]} \;\; &\( (6,W(s))\in B \equ (3,W(s))\notin B \).
\end{align*}
According to the instruction $b_7$, the set $B_7$ must encode the union of the sets encoded by $B_1$ and $B_2$:
\begin{equation*}
 \forall s\in\R^{[3]} \;\; \( (7,W(s))\in B \equ \( (1,W(s))\in B \lor (2,W(s))\in B \) \).
\end{equation*}
According to the instruction $b_8$, the set $B_8$ must encode the permutation of the set encoded by $B_1$:
\begin{equation*}
 \forall [x,y,z]\in\R^{[3]} \;\; \( (8,W([x,y,z]))\in B \equ (1,W([y,x,z]))\in B \);
\end{equation*}
similarly,
\begin{align*}
 \forall [x,y,z]\in\R^{[3]} \;\; &\( (9,W([x,y,z]))\in B \equ (1,W([x,z,y]))\in B  \); \\
 \forall [x,y,z]\in\R^{[3]} \;\; &\( (10,W([x,y,z]))\in B \equ (2,W([y,x,z]))\in B \); \\
 \forall [x,y,z]\in\R^{[3]} \;\; &\( (11,W([x,y,z]))\in B \equ  (2,W([x,z,y]))\in B \).
\end{align*}
According to the instruction $b_{12}$, the set $B_{12}$ must encode the Cartesian product (by $\R$) of the set
encoded by $B_1$:
\begin{equation*}
 \forall [x,y,z,u]\in\R^{[4]} \;\; \( (12,W([x,y,z,u]))\in B \equ (1,W([x,y,z]))\in B \);
\end{equation*}
similarly,
\begin{align*}
 \forall [x,y,z,u]\in\R^{[4]} \;\; &\( (13,W([x,y,z,u]))\in B \equ (2,W([x,y,z]))\in B
  \); \\
 \forall [x,y]\in\R^{[2]} \;\; &\( (14,W([x,y]))\in B \equ (3,W([x]))\in B \).
\end{align*}
According to the instruction $b_{15}$, the set $B_{15}$ must encode the projection of the set encoded by $B_1$:
\begin{equation*}
 \forall [x,y]\in\R^{[2]} \;\; \( (15,W([x,y]))\in B \equ \exists z\in\R \;\; (1,W([x,y,z]))\in B \);
\end{equation*}
similarly,
\begin{equation*}
 \forall [x,y]\in\R^{[2]} \;\; \( (16,W([x,y]))\in B \equ \exists z\in\R \;\; (2,W([x,y,z]))\in B \).
\end{equation*}

Toward the general formulation. We observe that the first two cases (complement and union) are unproblematic, while the
other three cases (permutation, set multiplication, and projection) need some additional effort. The informal
quantifiers like ``$\forall [x,y,z]$'' should be replaced with ``$\forall s$'', and the needed relations between finite
sequences should be generalized (and formalized).

\begin{exercise}
 The binary relation of truncation $\{([x_1,\dots,x_{n+1}],[x_1,\dots,x_n]) \mid n\in\N,
 x_1,\dots,x_{n+1}\in\R\}$ is definable. Prove it.  \emph{Hint:} use the evaluation function.
\end{exercise}

% \begin{sloppypar}
\begin{raggedright}
\begin{exercise}
 The ternary relation of appendment $\{(x,[x_1,\dots,x_n],[x_1,\dots,x_{n+1}]) \mid n\in\N,
 x_1,\dots,x_n,x\in\R\}$ is definable. Prove it.
\end{exercise}
\end{raggedright}
% \end{sloppypar}

Now the reader should be able to compose himself the general formulation. Also the additional condition that prevents
irrelevant points should be stipulated. We conclude that the set $B = \{1\}\times B_1 \cup \{2\}\times B_2 \cup \dots
\subset \R^2$ is definable. For each $n$ we check, whether the relation encoded by $B_n$ is of the form $\{u\}$ for
$u\in\R$ or not; if it is, we take $x_n=u$, otherwise $x_n=0$.  We get the definable sequence that contains all numbers
first-order definable in $(\R; +,\times,\N)$. The next step (explained in Section~\ref{sect4}), readily formalized (via
the function $D$ from Section~\ref{sect3}), provides a definable number not contained in this sequence.

\section{Fast-growing sequences}
\label{sect7}

Looking at decimal digits of two real numbers, for example,
\begin{align*}
  x &= 0.62831\,85307\,17958\,64769\,25286\,76655\,90057\,68394\,33879\,87502\,11641\dots
  \\ % \,94988\,91846\,15632\,81257\,2417\,97256\,06965\,06842\,34135\dots
  y &= 0.65465\,36707\,07977\,14379\,82924\,56246\,85835\,55692\,08082\,39542\,45575\dots
     % \,15320\,30341\,52669\,17935\,39584\,09434\,80222\,78477\,78618\dots
\end{align*}
can you see, which one is ``more definable''? Probably not. (Answer: $y=\sqrt{3/7}$ is
\href{http://en.wikipedia.org/wiki/Algebraic number}{algebraic}, therefore first-order
definable in $(\R;+,\times)$, while $x=\pi/5$ is not.) Surprisingly, a kind of visualization of definability is
possible in an interesting special case. The number
\begin{equation*}
\sum_{n=0}^\infty 10^{-3^n} =
0.\mathbf10\mathbf100\,000\mathbf10\,00000\,00000\,00000\,0\mathbf1000\,00000\,00000\,00000\dots
\end{equation*}
is \href{http://en.wikipedia.org/wiki/Transcendental number}{transcendental} (that is, not algebraic). Moreover, every
number of the form $\sum_{n=1}^\infty 10^{-k_n}$ with $ k_n \in \N $, $
\lim_{n\to\infty} \frac{k_{n+1}}{k_n} > 2 $ is transcendental, which follows from
\href{http://en.wikipedia.org/wiki/Roth's theorem}{Roth's theorem}.

\begin{exercise}
 If $(k_n)_{n=1}^\infty$ is a definable sequence of natural numbers, strictly increasing\index{strictly increasing} (that is, $ k_1<k_2<\dots $),
 then the number $ \sum_{n=1}^\infty 10^{-k_n} $ is definable. Prove it both in the framework of Section \ref{sect3}
 (first-order definability in $(\R; +,\times,\N)$) and the framework of Section \ref{sect5} (second-order
 definability). \emph{Hint:} $ \forall i\in\N \;\; \( D(i,x)=1 \equ \exists n\in\N \;\, i=k_n \) $,
 and $ \forall i\in\N \;\; D(i,x) \le 1 $.
\end{exercise}

\begin{exercise}\label{*}
  If a number $ \sum_{n=1}^\infty 10^{-k_n} $ with $ k_n \in \N $, $ k_1<k_2<\dots $, is definable, then the sequence
  $(k_n)_{n=1}^\infty$ is definable. Prove it in the framework of Section \ref{sect5} (second-order
  definability). \emph{Hint:} for every $n$, $k_{n+1}$ is the least $k$ such that $ k>k_n \land D(k,x)=1 $.
\end{exercise}
\enlargethispage{4pt}
Note that the sequence $(k_n)_{n=1}^\infty$ is defined by its property, which works only in the second-order
framework. The first-order framework requires an explicit relation between $n$ and $k=k_n$. Nevertheless, the claim of
Exercise \ref{*} holds also in the  framework of Section \ref{sect3}.\footnote{%
  % Remark about a proof.
  It is easy to obtain the sequence $(n_k)_{k=1}^\infty$ out of the sequence
  $(s_k)_{k=1}^\infty$ of sums $ s_k = \sum_{i=1}^k \al_i $ of the digits $ \al_i = D(i,x) $. The problem is that in the
  first-order framework we cannot define $(s_k)_{k=1}^\infty$ just by the property ``$ \forall k \;\;
  s_{k+1}=s_k+\al_{k+1} $''. Yet, this obstacle is surmountable; we can computably encode by natural numbers all tuples of
  natural numbers. (A similar trick was used in Section \ref{sect6}.)}

Thus, in order to get a first-order undefinable but second-order definable real number, it is sufficient to find a first-order
undefinable but second-order definable strictly increasing sequence of natural numbers. This can be made similarly to
Sections \ref{sect4}, \ref{sect6}, replacing Cantor's diagonal argument with the following fact:
\begin{compactitem}
  \item For every sequence of sequences (of numbers) there exists a strictly increasing sequence (of numbers) that
    overtakes\index{overtake} all the given sequences (of numbers).
\end{compactitem}
The proof is immediate: take $ y_n = n+\max_{i,j\in\{1,\dots,n\}} x_{i,j} $ where the number $ x_{i,j} $ is the $i$-th
element of the $j$-th given sequence; then clearly $ y_n > x_{n,m} $ whenever $ n \ge m $.

\begin{exercise}% \label{*}
  If the ternary relation $\{(i,j,x_{i,j})\mid i,j\in\N\}$ is definable, then the binary relation $\{(n,y_n)\mid
  n\in\N\}$ is definable. Prove it. \emph{Hint:} $ y=y_n \equ \( ( \exists i \;\, \exists j \;\; (i\le n\land
  j\le n\land y-n=x_{i,j})) \land ( \forall i \;\, \forall j \;\; ( i\le n\land j\le n \imp y-n \ge x_{i,j} ) )
  \) $.
\end{exercise}

Reusing the construction of Section \ref{sect6}, we enumerate all sequences of natural numbers,
definable in the framework of Section \ref{sect3}, by enumeration definable in the framework of Section \ref{sect5}, and
then overtake them all by a strictly increasing sequence of natural numbers, definable in the framework of
Section \ref{sect5}.

\begin{sloppypar}
To fully appreciate the incredible growth rate of this sequence, we note that it overtakes all
computable sequences, as well as an extremely fast-growing sequence $(M_N)_{N=1}^\infty$ mentioned in Introduction.
Recall $A_N$ and $A_{M,N}$ discussed there.
In the framework of Section \ref{sect3}, the ternary relation $\{(M,N,A_{M,N}) \mid M,N \in \N\}$, being recursively
enumerable (therefore Diophantine) is definable; and the binary relation $\{(N,A_N) \mid N \in \N\}$ is definable, since
$ a=A_N \equ \( ( \exists M\in\N \;\; a=A_{M,N} ) \land ( \forall M\in\N \;\; a\ge A_{M,N} ) \)
$. Defining $M_N$ as the least $M$ such that $A_{M,N}=A_N$ we observe that the sequence $(M_N)_{N=1}^\infty$ is
definable (since $ m=M_N \equ ( A_{m,N}=A_N \land A_{m-1,N}<A_N ) $). On the other hand, as noted in
Introduction, this sequence cannot be bounded from above by a computable sequence. 
\end{sloppypar}

More discussions of large numbers are available, see Scott Aaronson,\endnote{%
 \href{http://en.wikipedia.org/wiki/Scott Aaronson}{Aaronson, Scott} (1999).
 \href{https://www.scottaaronson.com/writings/bignumbers.pdf}{``Who can name the bigger number?''} \emph{Self-published.}}
John Baez\endnote{%
 \href{http://en.wikipedia.org/wiki/John C. Baez}{Baez, John} (2012).
 \href{https://johncarlosbaez.wordpress.com/2012/04/24/enormous-integers/}{``Enormous integers''}. \emph{Self-published.}}
 and references therein. A quote from Aaronson (pages 11--12):
\begin{quote}
  \emph{You defy him to name a bigger number without invoking Turing machines or some equivalent. And as he ponders this
  challenge, the power of the Turing machine concept dawns on him.}
\end{quote}
Definability could be mentioned here along with Turing machines.

\section{Definable but uncertain}
\label{sect8}

Two sets are called \href{http://en.wikipedia.org/wiki/equinumerous}{equinumerous}\index{equinumerous} if there exists a one-to-one
correspondence between them. In particular, two subsets $A,B$ of $\R$ are equinumerous if (and only if) $\exists f \in
\oP(\R^2) \;\; \( f \subset A\times B \land ( \forall x\in A \;\, \exists! y\in B \;\; (x,y)\in f ) \land (
\forall y\in B \;\, \exists! x\in A \;\; (x,y)\in f ) \) $. We see that the binary relation ``equinumerosity'' on $
\oP(\R) $ is second-order definable.

Some subsets of $\R$ are equinumerous to $\{1,\dots,n\}$ for some $n\in\N$ (these are finite sets). Others may be
equinumerous to $\N$ (these are called \index{countable}countable, or countably infinite), or $\R$ (these are called sets of cardinality
continuum),\index{cardinality continuum} or\dots\ what else? Can a set be more than countable but less than continuum?

This seemingly innocent question is one of the most famous in set theory,\endnote{%
  \emph{``Are there any sizes in between?'' This question (in the case of a negative answer, it is the
  \href{http://en.wikipedia.org/wiki/continuum_hypothesis}{continuum hypothesis}), is one of the most famous in set
  theory. Much as in the case of the parallel postulate, it was widely believed that the continuum hypothesis could
  simply be proven from \href{http://en.wikipedia.org/wiki/Zermelo-Fraenkel_set_theory}{ZFC}, and Cantor and many others
  devoted enormous time and effort to developing such a proof. It was not until much later that the combined efforts of
  G\"odel and Cohen established once and for all:\\
  \textsc{Theorem 3 (G\"odel, Cohen)} The continuum hypothesis is independent of ZFC.}\quad\hbox{}\hfill (From: J.~Reitz
  \cite[Section 3]{Re}.)}
the first among the \href{http://en.wikipedia.org/wiki/Hilbert's problems}{Hilbert's problems}.
The answer was expected to be ``no such sets'', which is the
\href{http://en.wikipedia.org/wiki/continuum_hypothesis}{continuum hypothesis}\index{continuum hypothesis}
(CH);\index{CH, the continuum hypothesis} Georg Cantor tried hard to prove
it, in vain; Kurt G\"odel proved in 1940 that CH cannot be disproved within
\href{http://en.wikipedia.org/wiki/Zermelo-Fraenkel set theory}{the axiomatic set theory called ZFC},\index{ZFC,
  axiomatic set theory} and hoped that new axioms will disprove it;\endnote{%
  \emph{Therefore one may on good reason suspect that the role of the continuum problem in set theory will be this, that
    it will finally lead to the discovery of new axioms which will make it possible to disprove Cantor's conjecture.}
  (From: K.~G\"odel \cite[the end]{Go}.)\\
  \emph{It was G\"odel who first suggested that perhaps ``strong axioms of infinity'' (large cardinals) could decide
    interesting set-theoretical statements independent over ZFC, such as CH. This hope proved largely unfounded for CH
    --- one can show that virtually all large cardinals defined so far do not affect the status of CH.} (From: R.~Honzik
  \cite[Abstract]{Ho}.)}
Paul Cohen proved in 1963 that CH cannot be proved within ZFC, and felt intuitively that it is \emph{obviously}
false.\endnote{%
  \emph{A point of view which the author feels may eventually come to be accepted is that CH is} obviously
  \emph{false.} [\dots] \emph{This point of view regards $C$ as an incredibly rich set given to us by one bold new
  axiom, which can never be approached by any piecemeal process of construction. Perhaps later generations will see the
  problem more clearly and express themselves more eloquently.} (From: P.~Cohen \cite[p.~151]{Co}.)\\
  \emph{If we really believe that the set-theoretic universe has to be built up piecemeal we surely cannot accept an
    axiom according to which enormous new sets (enormous because there is a jump in cardinality) simply
    nonconstructively appear.} (From: Nik Weaver \cite[p. 5]{We}.)}
Nowadays some experts hope to find ``the missing axiom'', others argue that this is hopeless.\endnote{%
  \emph{Many set theorists yearn for a definitive solution of the continuum problem, what I call a dream solution, one
    by which we settle the continuum hypothesis (CH) on the basis of a new fundamental principle of set theory, a
    missing axiom, widely regarded as true, which determines the truth value of CH.} [\dots] \emph{If achieved, a dream
    solution to the continuum problem would be remarkable, a cause for celebration.}\\
  \emph{In this article, however, I shall argue that a dream solution of CH has become impossible to achieve. Specifically,
  what I claim is that our extensive experience in the set-theoretic worlds in which CH is true and others in which CH
  is false prevents us from looking upon any statement settling CH as being obviously true.}
  (From: J.~Hamkins \cite{Ha2}.)\\
  Such a situation with the continuum problem raises doubts about the widely used set theory, especially the
  \href{http://en.wikipedia.org/wiki/axiom of powerset}{axiom of powerset}. Alternative, more
  \href{http://en.wikipedia.org/wiki/Constructive set theory}{constructive} approaches attract the attention of
  mathematicians and philosophers \cite{Ra}, \cite{Lin}, \cite{Ko}, \cite{Pa}. Maybe the collection of all subsets of an
  infinite set should be treated as a class rather than a set.}

A wonder: million published theorems\endnote{%
  2.3 mln articles in science and engineering are published in 2014, of them $2.6\%=0.06$ mln in mathematics. (See
    \href{https://www.nsf.gov/statistics/2018/nsf18300/}{National Center for Science and Engineering
    Statistics}. 2018.) Also, 0.3 mln articles were submitted to arXiv till now, of them 0.03 mln in 2014. (See
  \href{http://arxiv.org/help/stats/2017_by_area}{arXiv submission rate statistics}. 2017.) Thus, probably, 0.6 mln math
  articles are published for now. I guess, the average number of theorems per article is at least 2, which gives 1.2 mln
  published theorems. Some of them are notable, many are so-so, some are not new. And probably, hundreds or even
  thousands of them only pretend to be theorems, because of unnoticed errors in proofs. On the other hand, numerous
  lemmas formally are theorems. True, authors usually build proofs in the framework of the relevant branch of
  mathematics; but nearly all branches are embedded into ZFC. \emph{Today ZFC is the standard form of axiomatic set
    theory and as such is the most common foundation of mathematics.} (From
  Wikipedia:\href{https://en.wikipedia.org/wiki/Zermelo-Fraenkel_set_theory}{ZFC}.)}
in all branches of mathematics formally are deduced from the 9 axioms of ZFC; they answer, affirmatively or negatively,
million mathematical questions; some questions remain \href{http://en.wikipedia.org/wiki/open problem}{open}, waiting
for solutions in the ZFC framework; but the continuum hypothesis is an exception!\endnote{%
 Most notable exception, not the only exception. About 30 exceptions are available in
 Wikipedia:\href{http://en.wikipedia.org/wiki/List of statements independent of ZFC}{List of statements independent of
 ZFC}.}

Back to definability. Consider the set $Z$ of all subsets of $\R$ that are more than countable but less than
continuum. We do not know, whether $Z$ is empty or not, but anyway, we know that $Z$ is second-order definable. We
define a number $z$ by the following property:
\begin{equation*}
  (\, z=0 \;\land\; \exists A \;\; A\in Z \,) \;\;\lor\;\; (\, z=1 \;\land\; \forall A \;\; A\notin Z \,) \, .
\end{equation*}
That is, $z$ is $1$ if CH is true, and $0$ otherwise. This is a valid definition; $z$ is second-order definable; but we
cannot know, is it $0$ or $1$. Each one of the two equalities, $z=0$ and $z=1$, could be added (separately!) to the
axioms of ZFC without contradiction;\footnote{%
  Assuming, of course, that ZFC itself is consistent.}
according to the model theory, it means existence of two models of ZFC, one with $z=0$, the other with $z=1$. In this
sense, $z$ is model dependent.\index{model dependent}

Is $z$ computable? Yes, it is, just because $0$ and $1$ are computable numbers, and $z$ is one of these. You might feel
bothered, even outraged, but this is a valid argument. Compare it with the well-known proof that
an irrational elevated to an irrational power may be rational: $\({\sqrt2}\)^{\sqrt2}$ is either rational (which gives the
needed example), or irrational, in which case $\({\sqrt2}^{\sqrt2}\)^{\sqrt2}=\({\sqrt2}\)^2=2$ gives the needed
example.\endnote{%
 Wikipedia, % Section ``Examples'' in
 \href{https://en.wikipedia.org/wiki/Law of excluded middle\#Examples}{``Law of excluded middle''}; also % Section ``Properties'' in
 \href{https://en.wikipedia.org/wiki/Gelfond-Schneider constant\#Properties}{``Gelfond-Schneider constant''}.}
Seeing this, some retreat to
\href{https://en.wikipedia.org/wiki/Intuitionism}{intuitionism}, but almost all mathematics is
\href{https://en.wikipedia.org/wiki/Classical mathematics}{classical}, it accepts the
\href{https://en.wikipedia.org/wiki/Law of excluded middle}{law of excluded middle} and cannot arbitrarily disallow
it in some cases.

So, what is the algorithm for computing $z$? Surely the definition of an algorithm disallows such condition as ``if CH
holds, then'' within an algorithm. However, it cannot disallow a model-dependent algorithm $ A = (\text{if CH holds
  then } A_1 \text{ else } A_0) $, where $A_1$ is a (trivial) algorithm that computes the number 1, and $A_0$ computes
$0$. The conditioning ``if CH holds, then'' is allowed outside the algorithms (similarly, the conditioning ``if
$\({\sqrt2}\)^{\sqrt2}$ is rational, then'' is allowed outside the formulas). If you are unhappy with the affirmative
answer to the question ``is $z$ computable?'', ask a different question: ``is $z$ computable by a model-independent
algorithm?'' The answer is negative (see below).

On the other hand, definability of the number $z$ is established by a kind of ``generalized algorithm'' able to process
second-order objects (real numbers, relations between these, and relations between relations; recall the ``program''
$(b_k)_k$ in Section \ref{sect6}). This ``generalized algorithm'' is model independent, but its output is model
dependent.

In contrast, the number $\pi$ is model independent;\index{model independent} for every rational number $r$ one of the two inequalities $\pi<r$,
$\pi>r$ is provable in ZFC. The same applies to the numbers $ \sqrt{2}, \varphi, e $ discussed in Introduction, since
each of these numbers can be computed by a model independent algorithm. If a number is computable by a model-independent
algorithm, then this number is both model independent and computable.

What about Chaitin's constant $\Om$? It is limit computable by a model-independent algorithm. Also, it is first-order
definable (in $(\R; +,\times,\N)$), and the first-order framework disallows questions (such as CH) about
arbitrary sets of numbers, thus, one might hope that $\Om$ is model independent. But it is not!

Here we need one more fact about $\Om$. The sequence $(A_N)_{N=1}^\infty$ of its binary digits is not just uncomputable,
that is, the set $\{N\mid A_N=1\}$ is  not just non-recursive, but moreover, this set belongs to ``the most important class
of recursively enumerable sets which are not recursive'',\endnote{%
 From \href{https://www.encyclopediaofmath.org/index.php/Creative_set}{Encyclopedia of Mathematics:Creative set}.}
the so-called \href{http://en.wikipedia.org/wiki/creative and productive sets}{creative sets}, or equivalently, complete
recursively enumerable sets. Basically, it means that this sequence contains answers to all questions of the form ``does
the natural number $n$ belong to the recursively enumerable set $A$?'' And in particular(!), all questions of the form
``can the statement $S$ be deduced from the theory ZFC?'', since in ZFC (and many other formal theories as well) the set
of (numbers of) provable statements is recursively enumerable. Taking $S$ to be the negation of something provable (for
instance, $0\ne0$) we get the question ``is ZFC consistent?'' answered by one of the binary digits $A_N$ of $\Om$, whose
number $N_{\text{ZFC}}$ can be computed; if this $A_{N_{\text{ZFC}}}$ is $0$, then ZFC is consistent; if
$A_{N_{\text{ZFC}}}$ is $1$, then ZFC is inconsistent. However, by a famous G\"odel theorem, this question cannot be
answered by ZFC itself! Assuming that ZFC is consistent we have $A_{N_{\text{ZFC}}}=0$, but this truth is not provable
(nor refutable) in ZFC. (In fact, it is provable in \href{https://en.wikipedia.org/wiki/Large cardinal}{ZFC+large
  cardinal axiom}.) Therefore, in some models of ZFC we have $A_{N_{\text{ZFC}}}=0$, in others $A_{N_{\text{ZFC}}}=1$,
which shows that $\Om$ is model dependent. Moreover, there are versions of $\Om$ such that \emph{every} binary digit of
$\Om$ is model dependent \cite{So}, \cite{Ca}.

Yet the (first-order) case of $A_{N_{\text{ZFC}}}$ is less bothering that the (second-order) case of $z$, since we
still believe that $A_{N_{\text{ZFC}}}=0$. Adding the axiom ``$A_{N_{\text{ZFC}}}=1$'', that is, ``ZFC is inconsistent''
to ZFC we get a theory that is consistent\footnote{%
  Assuming, of course, that ZFC itself is consistent.}
but not \href{https://en.wikipedia.org/wiki/Omega-consistent}{$\om$-consistent}. This strange theory claims existence of
a proof of ``$0\ne0$'' in ZFC, of a finite length $N_{0\ne0}$, this length being a natural number. And nevertheless,
this theory claims that $N_{0\ne0}>1$, $N_{0\ne0}>2$, $N_{0\ne0}>3$, and so on, endlessly.\footnote{%
 Beware of the elusive distinction between two phrases, ``for each $n$ it claims $N_{0\ne0}>n$'' and ``it claims
 $\forall n \;\> N_{0\ne0}>n$''.}
Every model of this strange theory contains \href{https://en.wikipedia.org/wiki/Non-standard model of arithmetic}{more
  natural numbers} than the usual $1,2,3,\dots\,$ In mathematical logic we must carefully distinguish between two
concepts of a natural number, one belonging to a theory, the other to its
\href{https://en.wikipedia.org/wiki/Metamathematics}{metatheory}. In particular, when saying ``for every rational number
$r$ one of the two inequalities $\pi<r$, $\pi>r$ is provable in ZFC'' we should mean that $|r|$ is the ratio of two
metatheoretical natural numbers.

Using as binary digits an infinite sequence of \emph{independent} ``yes/no'' parameters of models of ZFC we get a model
dependent definable number $w$ whose possible values are \emph{all} real numbers. More exactly, the following holds in
the metatheory: for every real number $x$ there exists a model of ZFC\footnote{%
  Assuming, of course, that ZFC itself is consistent.}
whose natural numbers (and therefore rational numbers) are the same as in the metatheory, and for every
rational number $r$ the inequality $ w>r $ holds in the model if and only if $ x>r $. Is this possible in the first or
second order framework? I do not know. But in the third order framework this is possible, as suggested by the
\href{https://en.wikipedia.org/wiki/Generalized continuum hypothesis}{generalized continuum hypothesis}.\endnote{%
 \emph{In set theory, we have the phenomenon of the universal definition. This is a property $\phi(x)$, first-order
  expressible in the language of set theory, that necessarily holds of exactly one set, but which can in principle
  define any particular desired set that you like, if one should simply interpret the definition in the right
  set-theoretic universe. So $\phi(x)$ could be defining the set of real numbers $x=\R$ or the integers $x=\Z$ or the
  number $x=e^\pi$ or a certain group or a certain topological space or whatever set you would want it to be. For any
  mathematical object $a$, there is a set-theoretic universe in which $a$ is the unique object $x$ for which $\phi(x)$.}\\
 \textsc{Theorem.} \emph{Any particular real number r can become definable in a forcing extension of the universe.}\\
 (From: J.~Hamkins \href{http://jdh.hamkins.org/the-universal-definition/}{``The universal
   definition''}. Self-published. 2017. See also  \cite{HW}.)\\
 We adapt this idea. For notation used here, see
 Wikipedia:\href{https://en.wikipedia.org/wiki/Ordinal number}{Ordinal number}, in particular, Section ``Ordinals and
 cardinals''; $ \om=\om_0 $; and by the way, $ \N\cup\{0\} = \{0,1,2,\dots\} $. Also, following Kunen \cite{Ku}, the
 relation $ 2^\al=\be $ between two initial ordinals  $\al$ and $\be$ is interpreted as the same relation between the
 corresponding cardinals. Thus, $ 2^\om = \mathfrak c $ is the initial ordinal of the cardinality of the continuum.\\
 For every sequence $ (a_n)_{n=1}^\infty $ of binary digits $ a_n \in \{0,1\} $ the set $ \{
 (\om_n,\om_{\om+1+a_1+\dots+a_n}) \mid n\in\N\cup\{0\} \} \>\cup\> \{(\om_{\om+1},\om_{2\om})\} $  is a
 Easton index function \cite[Section 8.4, Def.~4.1]{Ku}; \href{https://en.wikipedia.org/wiki/Easton forcing}{Easton
   forcing} \cite[Section 8.4, Th.~4.7; also Corollary 4.8]{Ku} gives a model of ZFC  such that $
 2^{\om_n}=\om_{\om+1+a_1+\dots+a_n} $ for all $n\in\N\cup\{0\}$, and $ 2^{\om_{\om+1}} = \om_{2\om} $ (and the model
 has the same cardinals as the metatheory). Note that in this model $ \om_n < \mathfrak c $ (since $ \mathfrak c =
 2^{\om_0} = \om_{\om+1} $), and $ 2^{\om_n} < 2^{\mathfrak c} $ (since $ 2^{\mathfrak c} = \om_{2\om} = \om_{\om+\om}
 $).\\
 We consider the second-order definable set $ \B $ of all disjoint unions $ B = \{0\}\times B_0 \cup \{1\}\times B_1
 \cup \dots \subset \R^2 $ of sets $ B_0, B_1, B_2, \dots \subset \R $ such that $B_0$ is equinumerous to $\N$, and for
 each $n\in\N\cup\{0\}$,
 \begin{compactitem}
  \item $B_n$ is equinumerous to some subset of $B_{n+1}$;
  \item $B_n$ is not equinumerous to $B_{n+1}$;
  \item every subset of $B_{n+1}$ not equinumerous to $B_{n+1}$ is equinumerous to some subset of $B_n$.
 \end{compactitem}
 (It means that $B_n$ is of cardinality $\om_n$.)
 We note that equinumerosity of $\oP(B_n)$ and $\oP(B_{n-1})$ is third-order definable.
 If $\B$ is empty, we let $w=0$. Otherwise, for every $ n \in \N $ we
 define $b_n$ to be $0$ if $\oP(B_n)$ and $\oP(B_{n-1})$ are equinumerous for some (therefore,
 all) sets $ B = \{0\}\times  B_0 \cup \{1\}\times B_1 \cup \dots \in \B $; otherwise $b_n=1$.
 (It means that $b_n=a_n$.)
 Finally, we choose a definable map $ f $ from $[0,1]$ \emph{onto} $\R$ and let $ w = f(\sum_{n=1}^\infty 2^{-n} b_n)
 $.}

\section{Higher orders; set theory}
\label{sect9}

Recall the transition from first-order definability to second-order definability (Section \ref{sect5}); from the set $\R$
of all real numbers to the set $\mathbf S$ of all tuples and relations over $\R$, and the \Dstructure\ on $\mathbf S$
generated by the \Dstructure\ on $\R$ and two relations, membership and appendment, on $\mathbf S$. The next step
suggests itself: the set $ \mathbf T = \( \mathbf S \cup \mathbf S^2 \cup \mathbf S^3 \cup \dots \) \cup \(
\oP(\mathbf S) \cup \oP(\mathbf S^2) \cup \oP(\mathbf S^3) \cup \dots \) $ of all
tuples and relations over $\mathbf S$, with the \Dstructure\ on $\mathbf T$ generated by the \Dstructure\ on $\mathbf S$
and two relations, membership and appendment, on $\mathbf T$, formalizes third-order definability. This way we may
introduce infinitely many orders of definability, $ \R \subset \mathbf S \subset \mathbf T \subset \dots $, or $ T_1
\subset T_2 \subset T_3 \subset \dots $ where $ T_1 = \R $, $ T_2 = \mathbf S $, $ T_3 = \mathbf T $ and so
on. Similarly to Section \ref{sect6} we can prove that each order brings new definable real numbers (and new, faster-growing
sequences of natural numbers, recall Section \ref{sect7}).

But this is only the tip of the iceberg. The union of all these sets, $ T_\infty = T_1 \cup T_2 \cup \dots $, endowed
with the \Dstructure\ generated by the given \Dstructure s on all $T_n$, formalizes a new, transfinte order of
definability, and starts a new sequence of orders. Should we denote them by $T_{\infty+1}, T_{\infty+2}, \dots \,$? What
about $T_{\infty+\infty}\,$? How high is this hierarchy? Is it countable, or not?

Transfinite hierarchies are investigated by set theory (see Wikipedia:\href{https://en.wikipedia.org/wiki/Set
  theory}{Set theory}, and Section ``Some ontology'' there). Surprisingly, set theory does not need the field $\R$ of
real numbers as the starting point; not even the set $\N$ of natural numbers. A wonder: set theory is able to start from
nothing and get everything!\endnote{%
 \emph{So ecumenical set theorists instead spin this amazing structure from only the set that does not depend on the
 existence of anything: the empty set. This is the closest mathematicians get to creation from nothing!} (From
 \href{https://plato.stanford.edu/entries/nothingness/}{``Nothingness''}. Stanford Encyclopedia of philosophy. 2017.)}

\href{https://en.wikipedia.org/wiki/Von Neumann universe}{The cumulative hierarchy}\index{cumulative hierarchy} starts with the
\href{https://en.wikipedia.org/wiki/Empty set}{empty set}, denoted by $\emptyset$ or $\{\}$, the number $0$ defined as
just another name of the empty set, and stage zero, denoted by $V_0$ and defined as still another name of the empty
set. On the next step we consider the set $\oP(V_0)$ of all subsets of $V_0$. There is only one subset of
$\emptyset$, the empty set itself, thus $\oP(V_0) = \oP(\emptyset)=\{\emptyset\}=\{0\} $; we
define the number $1$ to be $\{0\}$, and stage one $ V_1=\oP(V_0) $. Similarly, $ \oP(V_1) $ is
the two-element set $ 2 = \{\emptyset,\{\emptyset\}\} = \{0,1\} = V_2 $, stage two. Somewhat dissimilarly, $
\oP(V_2) $ is the four-element set $ \{\emptyset,\{0\},\{1\},\{0,1\}\} = \{0,1,\{1\},2\} $, its
three-element subset $ \{0,1,2\} $ is (by definition) the number $3$, and $ V_3 = \oP(V_2) $ is the third
stage. More generally, $ n+1 = \{0,1,\dots,n\} \subset \oP(V_n) = V_{n+1} $ for $ n=1,2,3,4,\dots $. Thus,
$ V_n $ is a set of $\underbrace{2^{2^{\cdot^{\cdot^{2}}}}}_{n-1}$ elements(!), while $n$ is its subset of $n$
elements.

\begin{exercise}% \label{*}
  $ V_0 \subset V_1 \subset V_2 \subset \dots \, $ Prove it. \emph{Hint:} $ A \subset B \imp \oP(A) \subset \oP(B) $.
\end{exercise}

Here we face crossroads. One way is to treat the union $ V_0 \cup V_1 \cup V_2 \cup \dots $ of all $ V_n $ as the class
$V$ of all sets (a \href{https://en.wikipedia.org/wiki/Class (set theory)}{proper class}, not a set). This way leads to
the finite set theory (see Takahashi \cite{Tak}, Baratella and Ferro \cite{Bar} and others; see also
Wikipedia:\href{https://en.wikipedia.org/wiki/General set theory}{General set theory}). The other way is to treat the
union $ V_0 \cup V_1 \cup V_2 \cup \dots $ of all $ V_n $
as an infinite set, its infinite subset $\om = \{0,1,2,\dots\}$ as the first transfinite ordinal number, and $ V_\om =
\cup_{n=0}^\infty V_n $ as the first transfinite stage of the cumulative hierarchy. This way leads to the set theory
widely accepted by the mainstream mathematics.\endnote{%
  \emph{From the standpoint of mainstream mathematics, the great foundational debates of the early twentieth century
    were decisively settled in favor of Cantorian set theory, as formalized in the system ZFC (Zermelo-Fraenkel set
    theory including the axiom of choice). Although basic foundational questions have never entirely disappeared, it
    seems fair to say that they have retreated to the periphery of mathematical practice. Sporadic alternative proposals
    like topos theory or Errett Bishop's constructivism have never attracted a substantial mainstream following, and
    Cantor's universe is generally acknowledged as the arena in which modern mathematics takes place.} (From:
    Nik Weaver \cite[the first paragraph]{We}.)} 

\subsection{Finite set theory}
\label{sect9.1}

The finite set theory is equivalent (in some sense) to \href{https://en.wikipedia.org/wiki/Peano axioms}{arithmetic}
(Kaye and Wong \cite{Ka}); consistency of these theories is nearly indubitable,\endnote{%
  \emph{The vast majority of contemporary mathematicians believe that Peano's axioms are consistent, relying either on
    intuition or the acceptance of a consistency proof such as Gentzen's proof. A small number of philosophers and
    mathematicians, some of whom also advocate ultrafinitism, reject Peano's axioms because accepting the axioms amounts
    to accepting the infinite collection of natural numbers.} (From Wikipedia:\href{https://en.wikipedia.org/wiki/Peano
    axioms\#Consistency}{Peano axioms}.)}
in contrast to the (full) set theory whose \href{https://en.wikipedia.org/wiki/Axiom of infinity}{axiom of infinity} says
basically that the class $V_\om$ is a set.\endnote{%
  \emph{``There is something profoundly unsatisfactory about the axiom of infinity. It cannot be described as a truth of
    logic in any reasonable use of this term and so the introduction of it as a primitive proposition amounts in effect
    to the abandonment of Frege's project of exhibiting arithmetic as a development of logic'' (Kneale and Kneale,
    p.~699).\\
    This patched up set theory could not be identified with logic in the philosophical sense of ``rules for correct
    reasoning.'' You can build mathematics out of this reformed set theory, but it no longer passes as a foundation, in
    the sense of justifying the indubitability of mathematics.} (From \cite[p.~148--149]{He}.)}

In the finite set theory, the class $\N\cup\{0\}$ of numbers $0,1,2,\dots$ may be defined as the class of all sets $x$
such that $x$ is transitive, that is, $ \forall y \; \forall z \;\> ( y\in x \land z\in y \impl z\in x ) $, and $x$ is
totally ordered by membership, that is, $ \forall y \; \forall z \;\> \( ( y\in x \land z\in x ) \impl ( y=z \lor y\in z
\lor z\in y ) \) $. Adding the condition that $x$ is non-empty, that is, $ \exists y \; y \in x $, we get the class $\N$
of natural numbers $\{1,2,\dots\}$.

Each set $x$ is equinumerous to one and only one $n\in\N\cup\{0\}$; as before, ``equinumerous'' means existence of a set
$f$ such that $ f \subset x\times n $, that is, $ \forall p\in f \; \exists y\in x \; \exists m\in n \;\> p=(y,m) $
where $ (y,m)=\{\{y\},\{y,m\}\} $), and $f$ is a one-to-one correspondence between $x$ and $n$, that is, $ ( \forall
y\in x \;\, \exists! m\in n \;\; (y,m)\in f ) \land ( \forall m\in n \;\, \exists! y\in x \;\; (y,m)\in f ) $. In this
case we say that $n$ is the number of members of $x$.

The sum $m+n$ of $m,n\in\N\cup\{0\}$ may be defined as the number of members in the disjoint union $ \{1\}\times m \cup
\{2\}\times n $. 
The product $mn$ of $m,n\in\N\cup\{0\}$ may be defined as the number of members in the set product $ m \times n =
\{(k,\ell)\mid k\in m, \ell\in n\} $.
The power $m^n$ for $m,n\in\N\cup\{0\}$ may be defined as the number of functions from $n$ to $m$ (that is, from
$\{0,\dots,n-1\}$ to $\{0,\dots,m-1\}$).

A rational number could be defined as an equivalence class of triples $(p,n,q)$ of natural numbers $p,n,q\in\N$
w.r.t.\ such equivalence relation: $(p_1,n_1,q_1)\sim(p_2,n_2,q_2)$ when $ p_1 q_2 + n_2 q_1 = p_2 q_1 + n_1 q_2 $
(informally this means that $ \frac{p_1-n_1}{q_1} = \frac{p_2-n_2}{q_2} $, of course). However, in this case we cannot
introduce the class of rational numbers (since a proper class cannot be member of a class). Thus, it is better to choose
a single element in each equivalence class, and define a rational number as a triple $(p,n,q)$ of numbers
$p,n,q\in\N\cup\{0\}$ such that $q\ne0$, at least one of the two numbers $p,n$ is $0$, and the other is
\href{https://en.wikipedia.org/wiki/Coprime integers}{coprime} to $q$ (or $0$). We get the class $\Q$ of all rational numbers.
And, in order to treat natural numbers as a special case of rational numbers, we identify each natural number $n\in\N$
with the corresponding rational number $(n,0,1)\in\Q$.\endnote{%
 \emph{In many cases of interest there is a standard (or ``canonical'') embedding, like those of the natural numbers in
   the integers, the integers in the rational numbers, the rational numbers in the real numbers, and the real numbers in
   the complex numbers. In such cases it is common to identify the domain $X$ with its image $f(X)$ contained in $Y$, so
   that $X\subset Y$.} (From Wikipedia:\href{https://en.wikipedia.org/wiki/Embedding}{Embedding}.)}

Back to definability. We want to endow the class $V$ (of all sets in the finite set theory) with the \Dstructure\
generated by the membership relation $ \{ (x,y) \mid x\in y \} $. True, the notion of a \Dstructure\ on $V$ transcends
the finite set theory, since a collection of classes is neither a set nor a class. But still, in the metatheory, a class
may be called definable when it is obtainable from the membership relation by the 5 operations (complement, union,
permutation, Cartesian product, projection) introduced in Section \ref{sect2} for relations on the real line $\R$ in
particular, and arbitrary set in general. However, a pair of real numbers is not a real number, while a pair of sets is
a set! That is, $\R$ and $\R^2$ are disjoint; in contrast, $ V^2 \subset V $.
The order relation ``$x<y$'' between real numbers $x,y\in\R$ is a subset of $\R^2$ (rather than $\R$). But what about
the membership relation ``$x\in y$'' between sets $x,y\in V\,$? Should we treat it as a subclass of $V$ or $V^2\,$?

True, the plane $\R^2$ is not a subset of the line $\R$, but it can be injected into $\R$ by a definable function;
recall the injection $ W_2 : \R^2 \to \R $ introduced in the end of Section \ref{sect3} and used in Section \ref{sect6} for
encoding binary relations by unary relations. For example, the binary relation $ \{(x,y)\mid x<y\} $ is encoded by the
unary relation $\{W_2(x,y) \mid x<y\} $.

Here is a general lemma basically applicable in both situations, $W_2:\R^2\to\R$ and $ V^2\subset V$ (though, in the
latter case it needs some adaptation to the proper class).

\textbf{Lemma.}  Let $R$ be a set endowed with a \Dstructure, and $ f : R^2 \to R $ a definable injection. Then a binary
relation $ A \subset R^2 $ is definable if and only if the unary relation $ f(A) \subset R $ is definable.

\begin{exercise}% \label{*}
  Prove this lemma. \emph{Hint:} ``If'': $ (x,y)\in A \equ \exists z \;\> \( f(x,y)=z \land z\in f(A) \)
  $. \; ``Only if'': $ z\in f(A) \equ \exists x,y \;\> \( f(x,y)=z \land (x,y)\in A \) $.
\end{exercise}

\begin{exercise}% \label{*}
  All relations over $V$ mentioned above are definable classes. Prove it. \emph{Hint:}
  the equality relation ``$x=y$'' is ``$\forall z \;\> ( z\in x \equ z\in y ) $''; the relation
  ``$\{x\}=y$'' is ``$\forall z \;\> ( z=x \equ z\in y ) $''; the ternary relation ``$\{x,y\}=z$'' is
  ``$\forall u \;\> \( (u=x \lor u=y) \equ u\in z \) $''; the ternary relation ``$(x,y)=z$'' is
  ``$\exists u,v,w \;\> ( \{x\}=u \land \{x,y\}=v \land \{u,v\}=z )$''; the lemma applies; further, $ \N\cup\{0\} $ is
  the intersection of the class of transitive sets and the class of sets totally ordered by membership, etc. etc., up to
  ``$(p_1,n_1,q_1)\sim(p_2,n_2,q_2)$''.
\end{exercise}

Similarly, the basic relations between rational numbers are definable classes.

Real numbers cannot be represented by finite sets, but can be represented by classes (of finite sets) in several
ways. In the spirit of \href{https://en.wikipedia.org/wiki/Dedekind cut}{Dedekind cuts} we treat a real number as the
class of all rational numbers smaller than this real number. More formally: a real number is a subclass $A$ of $\Q$
such that
\begin{compactitem}
\item $A$ is a \href{https://en.wikipedia.org/wiki/Upper set}{lower class}; that is, $\forall a,b\in\Q \;\> ( a<b
  \land b\in A \impl a\in A ) $;
\item $A$ contains no greatest element; that is, $\forall a\in A \; \exists b\in A \;\> a<b $;
\item $A$ is not empty, and not the whole $\Q$; that is, $\exists a\in\Q \;\> a\in A $ and $\exists b\in\Q \;\>
  b\notin A $.
\end{compactitem}
And, in order to treat rational numbers as a special case of real numbers, we identify each rational number $a\in\Q$
with the corresponding real number $\{b\in\Q\mid b<a\} \in \R$.

Some examples. The real number $\sqrt2$ (``the Pythagoras' constant'') is the class of all rational numbers $a$ such that
$ a<0 \;\lor\; a^2<2 $.
The golden mean $\varphi$ is the class of all rational numbers $a$ such that $ a\le 0 \;\lor\; 0<a<1+\frac1a $.
The real number $e$ is the class of all rational numbers $a$ such that $ \exists n\in\N \;\>
\frac{(n+1)^n}{n^n} > a $.

Can we define $e$ via factorials, as in Exercise \ref{5.13}? We can define factorials without recursion; $n!$ is the
number of bijective functions from $n$ to itself (that is, from $\{0,\dots,n-1\}$ to itself; in other words,
permutations). But still, we need recursion when defining partial sums of the series $ \sum_{n=1}^\infty \frac1{n!} $
for $e$. Generally, an infinite sequence of rational numbers $ (s_n)_{n=1}^\infty $ is the class of pairs
$(n,s_n)$. Specifically, the sequence of partial sums of $\sum_n a_n$ is the class $S$ of pairs such that $ \forall
n\in\N \; \forall b\in\Q \;\> \( (n-1,b)\in A \impl (n,b+a_n)\in A \) $ (and $ \forall n\in\N \; \exists! b\in\Q \;\>
(n-1,b)\in A $, and $ (0,0)\in A $, of course). But we cannot define a class by its property! We deal with a
\Dstructure\ on $V$. A class must be defined by a common property of all its members, not a property of the class. Otherwise it would
be second-order definability in $V$ (thus, a transfinite level of the cumulative hierarchy). Can we formulate the
appropriate property of a pair $(n,s_n)$ alone? Yes, we can. Here is the property: there exists a function $ f :
\{0,\dots,n\} \to \Q $ such that $ f(0)=0 $ and $ \forall k\in\{1,\dots,n\} \;\> f(k)=f(k-1)+a_k $. The clue is that a
finite segment of the infinite sequence (of partial sums) is enough.

Similarly, an infinite sequence $(x_n)_{n=0}^\infty$ of sets $ x_n \in V $ is the class of pairs $(n,x_n)$, and it can
be defined recursively, by a recurrence relation of the form $ \forall n\in\N \;\> (x_{n-1},x_n) \in B $ where $ B $ is
a definable class of pairs, and an initial condition for $ x_0 $. (Use finite segments of the infinite sequence.)

Thus, every computable sequence of natural (or rational) numbers is a definable class of pairs. No need to use
Diophantine sets. Rather, for every \href{https://en.wikipedia.org/wiki/Turing machine}{Turing machine}, all possible
\href{https://en.wikipedia.org/wiki/Turing machine\#Turing_machine_"state"_diagrams}{``complete configurations''}
(called also ``situations'' and ``instantaneous descriptions'') may be treated as elements of a subclass of $V$, and the
rule of transition from one complete configuration to the next complete configuration may be treated as a definable
class of pairs (of complete configurations).

It follows that every computable real number, and moreover, every limit computable real number is definable. Having a
convergent definable sequence $(a_n)_n$ of rational numbers, we define its limit as the class of rational numbers $ b $
such that $ \exists n \; \forall k \;\> ( k>n \impl a_k>b+\frac1n ) $. In particular, $\pi$ (the Archimedes' constant)
and $\Om$ (the Chaitin's constant) are definable.

A sequence $(x_n)_n$ of real numbers cannot be treated as the class of pairs $(n,x_n)$ (since $x_n$ is not a set), but
can be treated as the disjoint union $ \{1\}\times x_1 \cup \{2\}\times x_2 \cup \dots $, that is, the set of pairs
$(n,a)$ where $a\in x_n$ (recall a similar workaround in Section \ref{sect6}). Also, a continuous function $ f : \R \to \R
$ cannot be treated as the class of pairs $\(x,f(x)\)$, but can be treated as the class of pairs $(a,b)$ of rational
numbers such that $ b < f(a) $. Such precautions allow us to translate basic calculus into the language of finite set
theory. However, arbitrary functions $ \R \to \R $ and arbitrary subsets of $\R$ are unavailable. Thus, the continuum
hypothesis makes no sense. Also, transferring
\href{https://en.wikipedia.org/wiki/Measure (mathematics)}{measure theory} and related topics (especially,
\href{https://en.wikipedia.org/wiki/Stochastic process}{theory of random processes}) to this ground (as far as possible)
requires effort and ingenuity.

The finite set theory can provide a reliable alternative airfield for much (maybe most) of the mathematical results
especially important for applications, in case of catastrophic developments in the transfinite hierarchy. Several
possible such ``alternative airfields'' are examined by mathematicians and philosophers \cite{Fe1}, \cite{Fe2},
\cite{We2}, \cite{Si2}, \cite{Si}, \cite{We3}, \cite{St}.

Informally, the finite set theory uses (for infinite classes) the idea of \href{https://en.wikipedia.org/wiki/Potential
  infinity}{potential infinity},\index{infinity!potential} prevalent before \href{https://en.wikipedia.org/wiki/Georg Cantor}{Georg Cantor}, while
the transfinite hierarchy uses the idea of \href{https://en.wikipedia.org/wiki/Actual infinity}{actual (completed)
  infinity},\index{infinity!actual} prevalent after Georg Cantor.\endnote{%
 \emph{In the philosophy of mathematics, the abstraction of actual infinity involves the acceptance (if the axiom of
   infinity is included) of infinite entities, such as the set of all natural numbers or an infinite sequence of
   rational numbers, as given, actual, completed objects. This is contrasted with potential infinity, in which a
   non-terminating process (such as "add 1 to the previous number") produces a sequence with no last element, and each
   individual result is finite and is achieved in a finite number of steps.} (From
  Wikipedia:\href{https://en.wikipedia.org/wiki/Actual infinity}{Actual infinity}.)\\
 \emph{Before Cantor, the notion of infinity was often taken as a useful abstraction which helped mathematicians reason
   about the finite world; for example the use of infinite limit cases in calculus. The infinite was deemed to have at
   most a potential existence, rather than an actual existence.} (From
  Wikipedia:\href{https://en.wikipedia.org/wiki/Controversy over Cantor's theory\#Reception_of_the_argument}{Controversy
    over Cantor's theory}.)\\
  \emph{The difficulty with finitism is to develop foundations of mathematics using finitist assumptions, that
   incorporates what everyone would reasonably regard as mathematics (for example, that includes real analysis).}
 (Ibid., \href{https://en.wikipedia.org/wiki/Controversy over Cantor's
    theory\#Objection_to_the_axiom_of_infinity}{``Objection to the axiom of infinity''}.)}

\subsection{Transfinite hierarchy}
\label{sect9.2}

The transfinite part of the cumulative hierarchy begins with \index{zzz@$\om$, the first transfinite ordinal}the first
transfinite ordinal number $ \om =
\{0,1,2,\dots\} $ (an infinite set) and the first transfinite stage $ V_\om = \cup_{n\in\om} V_n $ of the hierarchy (an
infinite set; $ \om \subset V_\om $). Note that $ x\in V_\om $ implies $ \oP(x) \in V_\om $, but $ x\subset
V_\om $ implies rather $ \oP(x) \subset V_{\om+1} $.
We continue as before:
\begin{align*}
  \om+1 &= \om \cup \{\om\} = \{0,1,2,\dots\} \cup \{\om\} \subset \oP(V_\om) = V_{\om+1} \, , \\
  \om+2 &= (\om+1) \cup \{\om+1\} = \{0,1,2,\dots\} \cup \{\om,\om+1\} \subset \oP(V_{\om+1}) = V_{\om+2}
\end{align*}
and so on; we get the stages $ V_{\om+n} $ for all finite $n$, and again, $ V_{\om+n} \subset V_{\om+n+1} $. The union
of all these stages is the stage $ V_{2\om} = V_\om \cup V_{\om+1} \cup V_{\om+2} \cup \dots = \cup_{\al<2\om} V_\al $
(still an infinite set), and $ 2\om = \{0,1,2,\dots\} \cup \{\om,\om+1,\om+2,\dots\} $ (an infinite subset of
$V_{2\om}$). Again, $ x\in V_{2\om} $ implies $ \oP(x) \in V_{2\om} $. Let us dwell here before climbing
higher.

Encoding of various mathematical objects by sets is somewhat arbitrary (see Wikipedia:
\href{https://en.wikipedia.org/wiki/Equivalent definitions of mathematical structures}{Equivalent definitions of
  mathematical structures}; likewise, an image may be encoded by files of \href{https://en.wikipedia.org/wiki/Image file
  formats}{type jpeg, gif, png etc.}), and their places in
the hierarchy vary accordingly. Treating a pair $(a,b)$ as $\{\{a\},\{a,b\}\}$ and a triple $(a,b,c)$ as $ \((a,b),c\) $
we get (for $0<n<\om$)
\begin{gather*}
  \forall a,b \in V_n \;\> (a,b) \in V_{n+2} \, ; \quad \forall a,b,c \in V_n \;\> (a,b,c) \in V_{n+4} \, ; \\
  \forall a,b \in V_\om \;\> (a,b) \in V_\om \, ; \quad \forall a,b,c \in V_\om \;\> (a,b,c) \in V_\om \, ; \\
  \forall a,b \in V_{\om+n} \;\> (a,b) \in V_{\om+n+2} \, ; \quad \forall a,b,c \in V_{\om+n} \;\> (a,b,c) \in
   V_{\om+n+4} \, .
\end{gather*}
Treating the set $\N$ of natural numbers as $\om\setminus\{0\}$ we get $ \N \subset V_\om $, $ \N \in V_{\om+1} $. 
Treating a rational number as an equivalence class of triples $(p,n,q)$ of natural numbers we get $ \Q \subset V_{\om+1}
$, $ \Q \in V_{\om+2} $, where $ \Q $ is the set of all rational numbers.
Alternatively, treating an integer as an equivalence class of pairs of natural numbers, and a rational number as an
equivalence class of pairs of integers,\endnote{%
  In Wikipedia see \href{https://en.wikipedia.org/wiki/Integer\#Construction}{Integer:Construction}, and
  \href{https://en.wikipedia.org/wiki/Rational number\#Formal_construction}{Rational number:Formal construction}.} 
we get
\begin{equation*}
  \Z \subset V_{\om+1} \, , \;\; \Z \in V_{\om+2} \, ; \quad
  \Q \subset V_{\om+4} \, , \;\; \Q \in V_{\om+5} \, ;
\end{equation*}
here $ \Z $\index{zzZ@$\Z$, the set of all integers} is the set of all integers. Treating a real number as a set of rational numbers we get
\begin{equation*}
  \R \subset V_{\om+n} \, , \quad \R \in V_{\om+n+1} \, ,
\end{equation*}
where $\R$ is the set of all real numbers, and $n$ is such that $ \Q \in V_{\om+n} $; be it $2$ or $5$, anyway, it
follows that $ \R \in V_{2\om} $.

Taking into account that generally $ A\in V_{2\om} \impl \oP(A) \in V_{2\om} $, and $ A,B\in V_{2\om} \impl
A\times B\in V_{2\om} $ (since $ A,B \subset V_{\om+n} \impl A\times B \subset V_{\om+n+2} $), we get $ \R^n \in
V_{2\om} $ and $ \oP(\R^n) \in V_{2\om} $ for all $n\in\N$. Every subset of $\R^n$ belongs to $V_{2\om}$,
and every set of subsets of $\R^n$ belongs to $V_{2\om}$; in particular, the
\href{https://en.wikipedia.org/wiki/Sigma-algebra}{\sia} of all \href{https://en.wikipedia.org/wiki/Lebesgue
  measure}{Lebesgue measurable} subsets of $\R^n$ belongs to $V_{2\om}$. Also, every function $ \R^n \to \R^m $ belongs
to $V_{2\om}$, and every set of such functions belongs to $V_{2\om}$; in particular, every equivalence class (under the
relation of equality \href{https://en.wikipedia.org/wiki/Almost everywhere}{almost everywhere}) of Lebesgue measurable
functions $ \R^n \to \R^m $ belongs to $V_{2\om}$, and the set \href{https://en.wikipedia.org/wiki/Lp
 space}{$L^1(\R^n)$} of all equivalence classes of \href{https://en.wikipedia.org/wiki/Lebesgue integration}{Lebesgue
 integrable} functions $ \R^n \to \R $ belongs to $V_{2\om}$. And the set of all
\href{https://en.wikipedia.org/wiki/Bounded operator}{bounded linear operators} $L^1(\R^n)\to L^1(\R^n)$ belongs to
$V_{2\om}$. Clearly, a lot of notable mathematical objects belong to $V_{2\om}$.\endnote{%
 \emph{$V_{\om+\om}$ is the universe of ``ordinary mathematics'', and is a model of Zermelo set theory.} (From
 Wikipedia:\href{https://en.wikipedia.org/wiki/Von Neumann universe\#Applications_and_interpretations}{Von Neumann
   universe}.)}

Would something like $V_{\om+100}$ suffice for all the objects mentioned above? The answer is negative as long as $ \R^n
$ is defined as $ \R^{n-1} \times \R = \{ (x,y) \mid x\in\R^n, y\in\R \} $ where $ (x,y) $ means $ \{\{x\},\{x,y\}\}
$. For every $ n\in\N $ the relation $ \R^n \notin V_{\om+2n-1} $ is ensured by the two exercises below.

\begin{exercise}% \label{*}
If $ A \times B \subset V_{\om+n+2} $, then $ A,B \subset V_{\om+n} $. Prove it. \emph{Hint:} $ \{\{a\},\{a,b\}\} =
(a,b) \in V_{\om+n+2} \impl a,b \in V_{\om+n} $.
\end{exercise}

\begin{exercise}% \label{*}
If $ A^{n+1} \subset V_{\om+2n} $ for some $n$, then $ A \subset V_\om $. Prove it. \emph{Hint:} induction in $n\ge1$,
and the previous exercise.
\end{exercise}

A more economical encoding is available (and was used in Section \ref{sect6}, see Exer.~\ref{6.1}, \ref{6.2}); instead of
the set $ \R^n $ of all $n$-tuples $(x_1,\dots,x_n)$ we may use the set $ \R^{[n]} $ of all $n$-sequences
$[x_1,\dots,x_n]$; as before, $ [x_1,\dots,x_n] = \{ (1,x_1), \dots, (n,x_n) \} $ is the set of pairs.

\begin{exercise}\label{9.6}
If $ A \in V_{\om+m+1} $, then $ A^{[n]} \in V_{\om+m+4} $ for all $ n\in\N $. Prove it. \emph{Hint:} $a_1,\dots,a_n \in
V_{\om+m} \impl [a_1,\dots,a_n] \in V_{\om+m+3} $.
\end{exercise}

A lot of theorems are published about real numbers, real-valued functions of real arguments, spaces of such functions
etc. I wonder, is there at least one such theorem sensitive to the distinction between $ V_{\om+100} $ and $
V_{\om+200} \,$? That is, theorem that can be formulated and proved within $ V_{\om+200} $ but not $ V_{\om+100} \,$? I
guess, the answer is negative. A seemingly similar question: is definability of real numbers sensitive to the
distinction between $ V_{\om+100} $ and $ V_{\om+200} \,$? I mean, is there at least one real number definable in $
V_{\om+200} $ but not $ V_{\om+100} \,$? This time, the answer is affirmative, as explained below.

For each $ n \in \N\cup\{0\} $ we endow the set $ V_{\om+n} $ with the \Dstructure\ $ D_{\om+n} $ generated by the
membership relation $ \{ (x,y) \mid x\in y \} $ (for $x,y \in V_{\om+n}$, of course).

Recall that, treating a real number as a set of rational numbers, and a rational number as an equivalence class of
triples $(p,n,q)$ of natural numbers, we have $ \Q \in V_{\om+2} $ and $ \R \in V_{\om+3} $. That is, $ \Q \subset
V_{\om+1} $ and $ \R \subset V_{\om+2} $.

Similarly to the finite set theory, $ \N $ and $ \Q $ are definable subsets of $ V_{\om+n} $ (whenever $n\ge1$), and the
basic relations between natural numbers are definable, as well as the basic relations between rational
numbers. Dissimilarly to the finite set theory, $ \R $ is a definable subset of $ V_{\om+n} $ (whenever $n\ge2$), and
the basic relations between real numbers are definable. An example: for $x,y\in\R $ we have $ x\le y \equ
\forall a\in\Q \;\> (a<x \impl a<y) \equ \forall a\in\Q \;\> (a\in x \impl a\in y) \equ x
\subset y $. Another example: for $x,y,z\in\R $ we have $ x+y=z \equ \forall c\in\Q \;\> \( c<z
\equ \exists a\in\Q \;\> ( a<x \land c-a<y) \) $.
Also the relation ``$x=\{b\in\Q\mid b<a\}$'' between a rational number $a$ and the corresponding real number $x$ is
definable, which implies definability of $\N$ embedded into $\R$. Thus, all real numbers first-order definable in
$(\R;+,\times,\N)$ (as in Section \ref{sect3}) are definable in $V_{\om+n}$ (whenever $n\ge2$).

What about second-order definability? It was treated in Section \ref{sect5} as a \Dstructure\ on the set $ \(
\cup_{n=1}^\infty \R^n ) \cup \( \cup_{n=1}^\infty \oP(\R^n) \) $, but a more economically encoded set $ S = \(
\cup_{n=1}^\infty \R^{[n]} ) \cup \( \cup_{n=1}^\infty \oP(\R^{[n]}) \) $ may be used equally well.

\begin{exercise}\label{9.7}
$ S \subset V_{\om+6} $. Prove it. \emph{Hint:} use Exercise \ref{9.6}.
\end{exercise}

Moreover, for every $n\ge6$, $S$ is a definable subset of $ V_{\om+n} $; and the four relations (that generate the
\Dstructure\ in Section \ref{sect5}) are definable relations on $ V_{\om+n} $. Thus, all real numbers second-order
definable as in Section \ref{sect5} are definable in $V_{\om+n}$ whenever $n\ge6$. In particular, the ``first-order
undefinable but second-order definable'' number of Section \ref{sect6} is definable in $V_{\om+6}$. However, all said does
not mean that it is undefinable in $V_{\om+2}$.

What we need is the second-order definability in $(V_{\om+2},D_{\om+2})$ rather than $(\R;+,\times,\N)$; that is,
definability in the set $ W_{\om+2} = \( \cup_{n=1}^\infty V^n_{\om+2} ) \cup \( \cup_{n=1}^\infty \oP(V^n_{\om+2}) \)
$.

\begin{exercise}% \label{9.7}
$ W_{\om+2} \subset V_{\om+6} $. Prove it. \emph{Hint:} similar to Exercise \ref{9.7}.
\end{exercise}

Once again, $ W_{\om+2} $ is a definable subset of $ V_{\om+6} $, and all real numbers definable in $ W_{\om+2} $ are
definable in $ V_{\om+6} $. That is, all real numbers second-order definable in $ V_{\om+2} $ are (first-order)
definable in $ V_{\om+6} $.

A straightforward generalization of Section \ref{sect6} gives a real number second-order definable in $ V_{\om+2} $ but
first-order undefinable in $ V_{\om+2} $. This number is definable in $ V_{\om+6} $ but undefinable in $ V_{\om+2}
$. Similarly, for each $ n\ge2 $ there exist real numbers definable in $ V_{\om+n+4} $ but undefinable in $ V_{\om+n}
$. We observe an infinite hierarchy of definability orders within $ V_{2\om} $.

Climbing higher on the cumulative hierarchy we get stages $ V_\al $ for ordinal numbers $ \al $ such as $ 2\om+n,
3\om+n, \dots $ Still higher, $ \om \cdot \om = \om^2 $, then $ \om^3, \dots $, then $ \om^\om, \om^{(\om^2)},
\om^{(\om^3)}, \dots \om^{(\om^\om)}, \dots $ Everyone may continue until feeling too dizzy; see
Wikipedia:\href{https://en.wikipedia.org/wiki/Ordinal notation}{Ordinal notation},
\href{https://en.wikipedia.org/wiki/Ordinal collapsing function}{Ordinal collapsing function},
\href{https://en.wikipedia.org/wiki/Large countable ordinal}{Large countable ordinal}. All these are countable
ordinals. By the way, every countable ordinal $\al$ may be visualized by a set of rational numbers, using a strictly
increasing function $ f : \al \to \Q $ (that is, $ f : \{\be\mid\be<\al\} \to \Q $). For example, $2\om$ may be
visualized by $ \{ 1-\frac1n \mid n\in\N \} \cup \{ 2-\frac1n \mid n\in\N \} $.

For every countable ordinal $\al\ge\om+2$ there exist real numbers definable in $ V_{\al+4} $ but undefinable in $ V_\al
$. Moreover, some of these real numbers are of the form $ \sum_{k=1}^\infty 10^{-k_n} $ (recall Section \ref{sect7}),
since there exists an increasing sequence (of natural numbers) definable in $ V_{\al+4} $ that overtakes all sequences
definable in $ V_\al $.

A wonder: stages $V_\al$ for $\al$ like $ \om^{\om^\om} $ are as far from ordinary mathematics as numbers like $
10^{10^{1000}} $ from ordinary engineering. Nevertheless these $ V_\al $ contribute to the supply of definable real
numbers.

Still higher, the set of all countable ordinals is the first uncountable ordinal $ \om_1 $. It cannot be visualized by a
set of rational or real numbers. Its cardinality is the first uncountable cardinality $ \aleph_1 $. The continuum
hypothesis is equivalent to the equality between $ \aleph_1 $ and the cardinality continuum.

For every ordinal $ \al \ge \om+2 $ (countable or not) the set of all real numbers definable in $ V_\al $ is countable
(and moreover, has an enumeration definable in $ V_{\al+4} $). In particular, the set of all real numbers definable in $
V_{\om_1} $ is countable. On the other hand, new definable real numbers emerge on all countable levels, and there are
uncountably many such levels. A contradiction?!

No, this is not a contradiction. Denoting by $R_\al$ the set of all real numbers definable in $V_\al$, and by
$\cO_\al$ the set of all ordinals definable in $V_\al$, we have $ R_\be \subset R_\al $ wherever $ \be \in \cO_\al $
(which follows from the lemma of Section \ref{sect5}). For all countable ordinals mentioned before we have $ \cO_\al = \al
$ (that is, all ordinals below $\al$ are definable in $V_\al$). In contrast, $ \cO_{\om_1} \ne \om_1 $, since $
\cO_{\om_1} $ is countable. The union $ \ti\R = \cup_{\al<\om_1} R_\al $ contains all real numbers definable \emph{with ordinal
  parameters} $ \al \in \om_1 $; but definability with parameters is outside the scope of this essay (recall
Section \ref{sect2}).

It is natural to ask, whether $ \cO_\al = \al $ for all countable ordinals $\al$, or not. Probably, we only know that
the affirmative answer cannot be proved without the choice axiom, and do not know, which answer (if any) can be proved
\emph{with} the choice axiom.\endnote{%
  In ZF (Zermelo-Frenkel set theory without the choice axiom) we can do the following.\\
  For arbitrary countable ordinal $\al$ the set $ D_\al $ of all $\al$-definable (that is, definable in $V_\al$)
  relations on $V_\al$ is countable, and
  has an $(\al+4)$-definable enumeration. Doing so for all $\al$ simultaneously we get a $\om_1$-definable function $ f
  : \om_1 \times \om \to \cup_{n=1}^\infty (\om_1)^n $ such that $ \forall \al \in \om_1 \;\> \forall A \in D_\al \;\>
  \exists n \in \om \;\; f(\al,n) = A $. Restricting ourselves to $A$ of the form $ \{\be\} $ (where $\be\in\om_1$) we
  get a $\om_1$-definable function $ g : \om_1 \times \om \to \om_1 $ such that $ \forall \al \in \om_1 \;\> \forall \be
  \in \cO_\al \;\> \exists n \in \om \;\; g(\al,n) = \be $.\\
  Now, for every sequence $ (\al_n)_{n\in\om} $ of countable ordinals, the union $ \cup_n \cO_{\al_n} $ is not just a
  countable union of countable sets, but a countable union of sets enumerated simultaneously by some function (which is
  trivial in ZFC but nontrivial in ZF), and therefore their union is countable, hence, not the whole $\om_1$. On the
  other hand, existence of $ (\al_n)_{n\in\om} $ such that $ \cup_n \al_n = \om_1 $ is consistent with ZF (Feferman and
  L\'evy 1963, see Cohen \cite[p.~143]{Co}). Thus, it is consistent with ZF that $ \cup_n \cO_{\al_n} \ne \cup_n \al_n
  $.\\
  This matter is closely related to a question about another definable (in $V_\al$, or otherwise) family, see
  MathStackExchange:\href{https://math.stackexchange.com/questions/2753698/is-it-possible-to-define-a-family-of-fundamental-sequences-for-all-countable-li}{``Is
    it possible to define a family of fundamental sequences for all countable (limit) ordinals? (Without AC)''},
  especially, the answer by Noah Schweber, and Remark 24 in Section 2.1 of Forster, Thomas E.
  \href{https://www.dpmms.cam.ac.uk/~tf/fundamentalsequence.pdf}{``A tutorial on countable ordinals''}. Self-published.\\
  Another question of this kind,
  \href{https://mathoverflow.net/questions/309494/is-the-smallest-l-alpha-with-undefinable-ordinals-always-countable}{``Is
    the smallest $L_\al$ with undefinable ordinals always countable?''} is answered affirmatively by
  \href{https://mhabic.github.io}{Miha Habi\v c.}}

If $\cO_\al = \al $ for all countable ordinals $\al$, then $R_\al \uparrow \ti R$, and $\ti R$
is uncountable (of the cardinality $\aleph_1$ of $\om_1$). Otherwise, there exists a countable ordinal $\al$ such
that $ \cO_\al \ne \al $ while $ \cO_\be = \be $ for all $\be < \al$ (and therefore $
\cO_\al = \be $ for some $\be < \al$). In this case the transfinite sequence $ (\cO_\al)_{\al<\om_1} $ is not monotone,
and we do not know, whether the union $\ti R$ is countable, or not.

Countability or uncountability of $\ti R$ matters for model dependence. There is a countable set of formulas (in the
language of the set theory) that define real numbers on all levels $V_\al$. Some are model independent, others are model
dependent. If $\ti R$ is countable, then each of these model dependent definable real numbers has at most countably many
possible values. Otherwise, if $\ti R$ is uncountable, then at least one of these model dependent definable real numbers
has uncountably many possible values.

This matter is closely related to the position of Laureano Luna \cite{Lu0} (see also \cite[pages 19--20]{Lu}):
\begin{quote}
  \emph{``Pieces of language taken as mere} syntactical \emph{expressions (letter-strings) should be
    distinguished from definitions, which are} semantical \emph{objects, namely,} interpreted \emph{letter-strings.''}
  (Page 61.)\\
\emph{``The meanings of the letter-strings that express definitions of reals are context-dependent, the context being
  here the definability level on which they are used. [\dots] if some letter-strings express more than one definition of
  a real number, there is no reason to think there are only countably many such definitions and only countably many
  definable real numbers.''} (Page 64.)
\end{quote}

Two objections arise. First, we did not prove that $\ti R$ is uncountable. Second, model dependence does not apply to
$V_\al$, since $V_\alpha$ is not a model of ZFC. We'll return to the second objection after climbing on the cumulative
hierarchy to $V_{\om_1}$ and much higher.

The stage $V_{\om_1}$ of the cumulative hierarchy is vast; its cardinality is very large (much larger than the
cardinality $\aleph_1$ of $\om_1$). Now consider the first ordinal $\al$ of this very large cardinality and the
corresponding stage $V_\al$. Iterating this jump we get a slight idea of the class of all sets, the incredible universe
$V_{\text{ZFC}}$ of the set theory ZFC. The whole $V_{\text{ZFC}}$ grows from a small seed, the first infinite ordinal
$\om$, whose existence is just postulated (the axiom of infinity).

If you want to soar above $V_{\text{ZFC}}$, you need a new axiom of infinity that ensures existence of
an ordinal $\al$ such that $V_\al$ is a model of ZFC; every such ordinal, being \href{https://en.wikipedia.org/wiki/Initial
  ordinal}{initial}, is a cardinal, called a \href{https://en.wikipedia.org/wiki/Worldly
  cardinal}{worldly cardinal}.\index{worldly cardinal} For climbing still higher try the so-called \href{https://en.wikipedia.org/wiki/Large
  cardinal}{large cardinals}. And be assured that these supernal stages do contribute to the supply of definable real
numbers.\endnote{%
  \emph{Eventually, we will in our definitions be attracted to the possibilities of using higher order mathematical
    objects and constructions, such as function classes, spaces or measures, and this amounts to defining objects in
    increasingly large fragments $V_\al$ of the set-theoretic universe. Most all of the classical mathematical structure
    is itself definable in the set-theoretic structure $\langle V_{\om+\om},\in\rangle$, a model of the Zermelo axioms,
    and so the definable reals of this structure includes almost every real ever defined classically. The structures
    arising with larger ordinals, however, allow us to define even more reals.} (From \cite[Section 1]{HLR}.)}

Assuming existence of large cardinals we get a transfinite hierarchy of worldly cardinals $\kappa_\al$, and
the corresponding models $W_\al = V_{\kappa_\al}$ of ZFC, for all countable ordinals $\al$ (and more; but here we do not
need uncountable $\al$). Using $W_\al$ instead of $V_\al$ we get new versions of $\cO_\al$, $R_\al$ and
$\ti R$. Again, we do not know, whether this new $\ti R$ is countable, or not. If it is uncountable, then again, at
least one model dependent definable real number has uncountably many possible values; and this time, model dependence
applies.

\subsection{Getting rid of undefinable numbers}
\label{sect9.3}

Climbing down to earth, is it possible to restrict ourselves to definable numbers and still use the existing theory of
real numbers and related objects? An affirmative answer was found in 1952 \cite{My} and enhanced recently \cite{HLR}.

Before climbing down we need to climb up to the first worldly cardinal $\al$ and the corresponding model $V_\al$ of ZFC.
Within the model we consider the \href{https://en.wikipedia.org/wiki/Constructible universe}{constructible hierarchy}
$(L_\be)_{\be\le\al}$, take the least $\be$ such that $ L_\be $ is a model of ZFC, and get the so-called
\href{https://en.wikipedia.org/wiki/Minimal model (set theory)}{minimal transitive model} of ZFC. This model is
pointwise definable\index{pointwise definable} \cite[``Minimal transitive model'']{Ha1}, that is, every member of this model is definable (in this
model).

Accordingly, this model is countable (and $\be$ is countable). Nevertheless, every theorem of ZFC holds in every model
of ZFC; in particular, \href{https://en.wikipedia.org/wiki/Georg Cantor's first set theory article}{Cantor's theorem}
``$\R$ is uncountable'' holds in the countable model $L_\be$. No contradiction; enumerations of $ \R \cap L_\be $ exist,
but do not belong to $L_\be$. Likewise, \href{https://en.wikipedia.org/wiki/Lebesgue measure}{a well-known theorem of
  measure theory} states that the interval $(0,1)$ cannot be covered by a sequence of intervals $(a_n,b_n)$ of total
length $\sum_{n=1}^\infty (b_n-a_n) < 1 $. True, for every $\eps>0$ the set $(0,1)\cap L_\be$, being countable, can be
covered by a sequence of intervals of total length $\eps$; but such sequences do not belong to $L_\be$ (even if
endpoints $a_n,b_n$ do belong). Working in $L_\be$ we have to ensure that all relevant objects (not only real numbers)
belong to $L_\be$.

\begin{quote}
$\bullet$ \emph{One often hears it said that since there are indenumerably many sets and only denumerably many names,
therefore there must be nameless sets. The above shows this argument to be fallacious.} (Myhill 1952, see \cite[the
  last paragraph]{My}.)
\end{quote}

\begin{quote}
$\bullet$ \emph{In my opinion, an object is conceivable only if it can be defined with a finite number of words.}
  (Poincar\'e 1910, translated from German, see \cite[page 58]{Lu0}.)
\end{quote}

\section{Conclusion}
\label{sect10}

Each definition (of a real number, or another mathematical object) is a finite text in a language. The language may be
formal (mathematical) or informal (natural). In both cases the text is composed of expressions that
\href{http://en.wikipedia.org/wiki/Reference}{refer} to objects and relations between objects. The
\href{http://en.wikipedia.org/wiki/Extension (semantics)}{extension} of an expression is the corresponding set of
objects, or set of pairs (of objects), or triples, and so on. For a mathematical language, all objects are
mathematical; a natural language may mention non-mathematical objects, and even itself, as in the phrase \emph{``The
  preceding two paragraphs are an expression in English that unambiguously defines a real number $r$''} (recall
Introduction, Richard's paradox), which leads to a problem: the mentioned expression in English fails to
define! \emph{``So when we speak in English about English, the `English' in the metalanguage is not exactly the same as
the `English' in the object-language.''} \cite[p.~15]{Lu}. A natural language, intended to be its own metalanguage, is
burdened with paradoxes. A mathematical language avoids such (and hopefully, any) paradoxes at the expense of being
different from its metalanguage. The metalanguage is able to enumerate all real numbers definable in the language and
define more real numbers.

On one hand, a natural language itself is inappropriate for mathematics. On the other hand, a mother tongue is always a
natural language. In order to avoid both restrictions of a fixed mathematical theory and paradoxes of a natural language
we may get the best of both worlds by considering two-part texts. The first part, written in a natural language,
introduces a mathematical theory. The second part, written in the (mathematical) language of this theory, defines some
real number.

In this framework the question ``is every real number definable?'' falls out of mathematics, because the notion
``mathematical theory'' above cannot be formalized. Some may admit only potential infinity and stop on the finite set
theory. Some may admit actual infinity and the transfinite hierarchy up to (exclusively) some preferred ordinal
(sometimes $2\omega$ \cite{Tar2}; more often, something controversially believed to be the first undefinable
ordinal\endnote{%
\emph{Thus, our definability levels cannot go beyond the countable definable ordinals. What these are is
  contentious. Constructivists will typically argue that all ordinals are constructive and that the (classical)
  \href{http://en.wikipedia.org/wiki/Church-Kleene ordinal}{least nonconstructive ordinal} $\om_1^{\text{CK}}$ does not
  exist, though from more liberal standpoints it is a definable countable ordinal. Some predicativists believe that the
  classical \href{http://en.wikipedia.org/wiki/Feferman–Schütte ordinal}{Feferman-Sch\"utte ordinal} $\Ga_0$ does not
  exist (though from the classical viewpoint it is a definable countable ordinal); some, as} Weaver 2009, \emph{believe
  predicativity can go beyond $\Ga_0$.} (Luna \cite[p. 64]{Lu0}.)}%
).
Or the whole universe of ZFC but no
more (equivalently, up to the first worldly cardinal).\endnote{%
 Someone may say: I just use ZFC ``as is'' and do not care about worldly cardinals. But in this case many
 interesting definable real numbers are model dependent (recall Section \ref{sect8}), they exploit large cardinals
 whenever possible.}
Or the \href{http://en.wikipedia.org/wiki/Tarski-Grothendieck set theory}{Tarski-Grothendieck set theory}. Or higher, up
to some preferred \href{http://en.wikipedia.org/wiki/Large cardinal}{large cardinal}. Or some more exotic
\href{http://en.wikipedia.org/wiki/Alternative set theory}{alternative set theory}. Or something brand new, like a kind
of \href{http://en.wikipedia.org/wiki/Homotopy type theory}{homotopy type theory}. Or even something not yet
published. Most choices mentioned above were unthinkable in the first half of the 20th century. Who knows what may
happen near year 2100? \emph{``Mathematics has no generally accepted definition''} (from
Wikipedia:\href{http://en.wikipedia.org/wiki/Definitions of mathematics}{Definitions of mathematics}); the same can be
said about ``mathematical theory''.\endnote{%
 The set of all \href{http://en.wikipedia.org/wiki/Godel's incompleteness theorems\#Formal_systems:_completeness,_consistency,_and_effective_axiomatization}{consistent effectively axiomatized
 formal theories} is well-defined but irrelevant, because most of them have no \href{http://en.wikipedia.org/wiki/Interpretation
 (logic)\#Intended_interpretations}{intended interpretation}.}
Admitting actual infinity we do not get rid of potential infinity; the latter returns, hardened, on a higher level
\cite{Si}, \cite{Li}. Maybe the recent, \href{http://logic.harvard.edu/EFI_Hamkins_Comments.pdf}{hotly debated}
conception of \href{http://en.wikipedia.org/wiki/Multiverse}{multiverse} is another fathomable segment of the
unfathomable potential infinity of mathematics. ``The most general definition of a definition'' appears to be as
problematic as ``the set of all sets''.

Another problem manifests itself as model dependence for a mathematical language, and context dependence for a natural
language. In both cases a single text may refer to different objects (depending on the model or the context,
respectively), which blurs the idea of definability. Recall Section \ref{sect8} (the last paragraph): every real number
is ``hardwired'' in some model of the set theory (ZFC). We may feel that definability of the number does not follow
unless the model is definable; but what do we mean by definability of a model? Another case, recall Section
\ref{sect9.2} (the last paragraph): uncountable hierarchy of models indexed by countable ordinals leads to a set of real
numbers, possibly uncountable (but maybe not). Nothing is ``hardwired'' here, except for these countable
ordinals. Outside mathematics, a more or less similar case, treated by Luna \cite{Lu0}, shows context dependence in a
hierarchy of contexts (levels of definability) indexed by definable countable ordinals. A mathematical counterpart with
uncountably many definable real numbers could exist in some alternative (to ZFC and alike; probably closer to
\href{http://en.wikipedia.org/wiki/Kripke–Platek set theory}{Kripke-Platek}) set theory such that the class of all
definable countable ordinals is not a countable set, and preferably, not a set at all \cite[p. 65]{Lu0}.

Bad news: definability is a very subtle property of a real number. Good news: other properties, more relevant to
applications, are unsubtle; and definability is rather of philosophical interest. \emph{``Mathematicians, in general, do
not like to deal with the notion of definability; their attitude toward this notion is one of distrust and reserve.''}
(Tarski \cite{Tar2}, the first phrase; now partially obsolete, partially actual.)

\setlength{\parskip}{1ex}
\def\enoteformat{\rightskip=0pt \leftskip=0pt \parindent=0em \leavevmode\llap{\makeenmark\;}}
\def\enotesize{\small}

\printendnotes

% \bigskip
% \newpage

\vfill

\printindex

\vfill

\end{document}

%% file: preamble.tex
\theoremstyle{definition}

\newtheorem{exercise}{Exercise}

\numberwithin{exercise}{section}
\numberwithin{equation}{section}

\renewcommand{\phi}{\varphi}

\newcommand{\ti}{\tilde}

\renewcommand{\(}{\bigl(}
\renewcommand{\)}{\bigr)\vphantom{)}}

\newcommand{\equ}{\Longleftrightarrow}

\newcommand{\impl}{\>\Longrightarrow\>}
\newcommand{\imp}{\Longrightarrow}

\newcommand{\oP}{\operatorname{P}}

\newcommand{\ccirc}{\mathbin{\kern-0.13em\circ}}

\newcommand{\eps}{\varepsilon}

\newcommand{\Ga}{\Gamma}
\newcommand{\om}{\omega}
\newcommand{\Om}{\Omega}

\newcommand{\al}{\alpha}
\newcommand{\be}{\beta}
\newcommand{\cO}{\mathcal O}

\newcommand{\B}{\mathcal B}

\newcommand{\R}{\mathbb R}

\newcommand{\Z}{\mathbb Z}
\newcommand{\N}{\mathbb N}
\newcommand{\Q}{\mathbb Q}

\newcommand{\Dstructure}{D\nobreakdash-\hspace{0pt}structure}
\newcommand{\sia}{$\sigma$\nobreakdash-\hspace{0pt}algebra}

\newcommand{\nbdh}[2]{$#1$\nobreakdash-\hspace{0pt}#2}

\newcommand{\ary}[1]{$#1$\nobreakdash-\hspace{0pt}ary}

\makeatletter
    \def\thebibliography#1{\section*{Bibliography\@mkboth
      {FOOTNOTES}{FOOTNOTES}}\list
      {[\arabic{enumi}]}{\settowidth\labelwidth{[#1]}\leftmargin\labelwidth
	\advance\leftmargin\labelsep
	\usecounter{enumi}}
	\def\newblock{\hskip .11em plus .33em minus .07em}
	\sloppy\clubpenalty4000\widowpenalty4000
	\sfcode`\.=1000\relax}
    \makeatother